\documentclass[review]{elsarticle}

\usepackage{hyperref}

\journal{Computers and Electronics in Agriculture}

\biboptions{authoryear}

\usepackage{amsmath}
\usepackage{amssymb}
\usepackage{epstopdf}
\usepackage{graphicx}
\usepackage{eurosym}
\usepackage{color}
\newtheorem{remark}{Remark}[section]

\begin{document}

\begin{frontmatter}

\title{A model-free control strategy \\ for an experimental greenhouse \\ with an application to fault accommodation}


\author[mymainaddress,mysecondaryaddress]{Fr\'ed\'eric Lafont\corref{mycorrespondingauthor}}
\cortext[mycorrespondingauthor]{Corresponding author.}
\ead[url]{http://www.lsis.org}

\author[mymainaddress,mysecondaryaddress]{Jean-Fran\c cois Balmat}

\author[mymainaddress,mysecondaryaddress]{\\ Nathalie Pessel}

\author[mythirdaddress,myfourthaddress]{Michel Fliess}

\address[mymainaddress]{Universit\'e de Toulon, CNRS, LSIS (UMR 7296), 83957 La Garde, France. \newline {\tt\small \{lafont, balmat, nathalie.pessel\}@univ-tln.fr}}
\address[mysecondaryaddress]{Aix Marseille Universit\'e, CNRS, ENSAM, LSIS (UMR 7296), 13397 Marseille, France.}
\address[mythirdaddress]{LIX (UMR CNRS 7161), \'Ecole polytechnique, 91128 Palaiseau, France.  \newline
        {\tt\small Michel.Fliess@polytechnique.edu}}
\address[myfourthaddress]{AL.I.E.N. (Alg\`{e}bre pour Identification \& Estimation Num\'{e}riques), \\ 24-30 rue Lionnois, BP 60120, 54003 Nancy, France. \newline
        {\tt\small michel.fliess@alien-sas.com}}

\begin{abstract}
Writing down mathematical models of agricultural greenhouses and regulating them via advanced controllers are challenging tasks since strong perturbations, like meteorological variations, have to be taken into account. This is why we are developing here a new model-free control approach and the corresponding ``intelligent'' controllers, where the need of a ``good'' model disappears. This setting, which has been introduced quite recently and is easy to implement, is already successful in many engineering domains. Tests on a concrete greenhouse and comparisons with Boolean controllers are reported. They not only demonstrate an excellent climate control, where the reference may be modified in a straightforward way, but also an efficient fault accommodation with respect to the actuators.
\end{abstract}


\end{frontmatter}

\section{Introduction}
\begin{figure}[!t]
      \centering
      \includegraphics[scale=0.5]{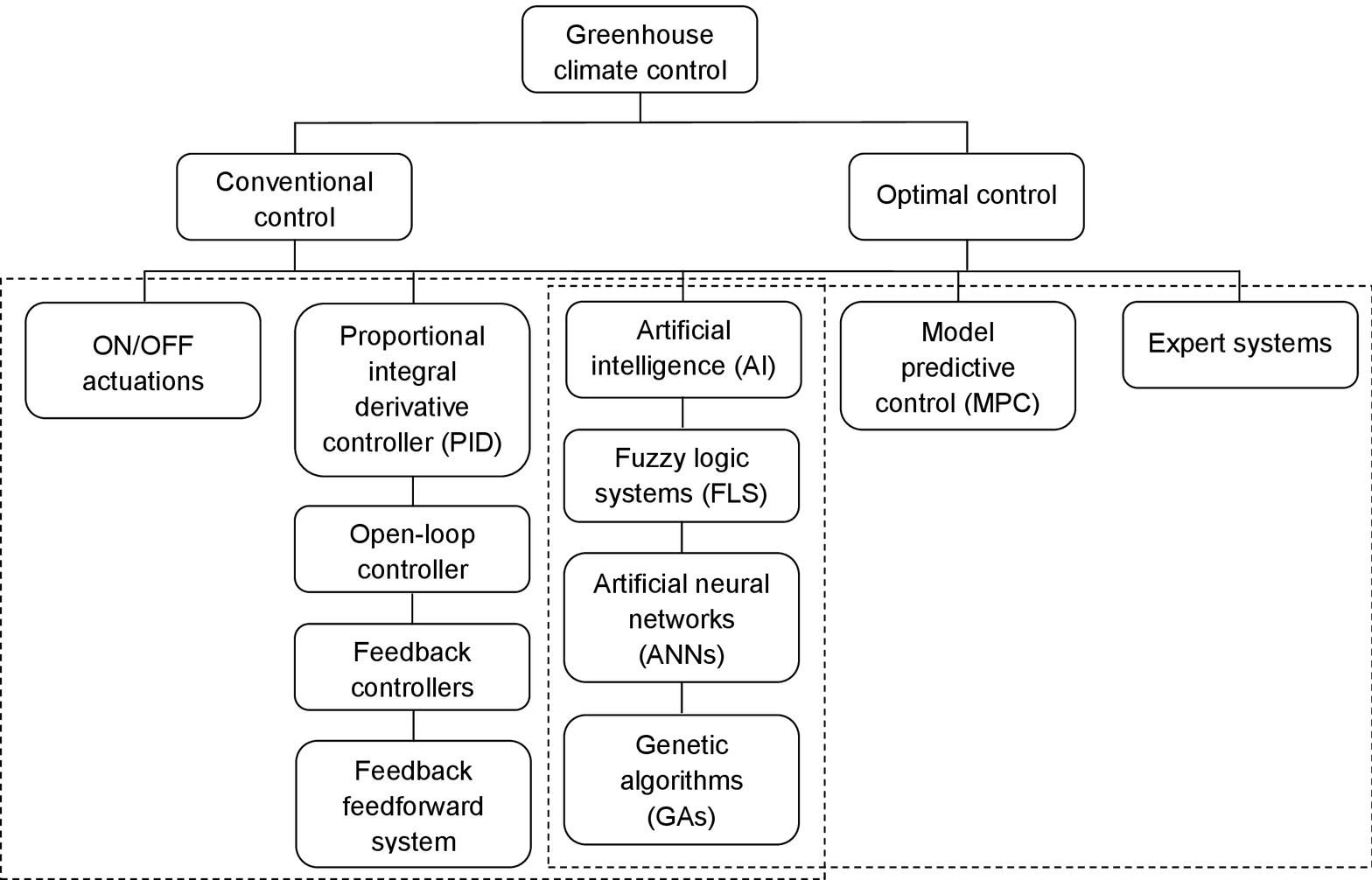}
      \caption{Greenhouse control theories classification in \cite{c19}}
      \label{Theories}
   \end{figure}

Table \ref{bilan} in \cite{Callais} shows that already a few years ago a large percentage of agricultural greenhouses were computerized.
\begin{table}[!tbp]
\caption{Percentage distribution of surfaces for the soilless crop greenhouses in France in 2005}
\label{bilan}
\begin{center}
\begin{tabular}{|c||c||c||c|}
\hline
 \multicolumn{4}{|c|}{Climate control}\\
\hline
Without & Manual & Automated & Computerized\\
\hline
6 \% & 7 \% & 20 \% & 67 \% \\
\hline
\end{tabular}
\end{center}
\end{table}
The corresponding automated microclimate regulation
should not only improve the production and its quality but also reduce pollution and energy consumption. Most of the existing control approaches, like adaptive control, predictive control, optimal control, stochastic control, nonlinear control, infinite dimensional systems, PIDs, On/Off, or Boolean, control, fuzzy control, neural networks, soft computing, expert systems, \dots, have been employed and tested. The literature on the modeling and control of greenhouses is therefore huge.
See, \emph{e.g.},:
\begin{itemize}
\item the books by \cite{eue,crc14,rod,tf,c15,von}; and the references therein,
\item the papers and memoirs by \cite{aa,c5,c8,c7,blasco,cap,udink,critten,c1,dong,c19,el,fourati,gruber,io,c16,kit,c6,pas,c18,c4,Pin,c2,sham,speetens,tcham,c17,zh}; and the references therein.
\end{itemize}
Let us summarize, perhaps too briefly, some of the various control aspects which were developed in the above references (see, also, Figure \ref{Theories}):
\begin{itemize}
\item writing down a ``good'' model, which is
necessarily nonlinear, either via physical laws or via black box identification, leads to most severe calibration and robustness issues, especially with respect to strong weather disturbances, which
are impossible to forecast precisely,
\item for multi-models appropriate control laws are difficult to synthesize,
\item ``conventional'' PID and On/Off techniques, {which preclude any mathematical modeling}, are therefore the most popular in industrial greenhouses, although:
\begin{itemize}
\item they are difficult to tune,
\item their performances are far from being entirely satisfactory.
\end{itemize}
\end{itemize}

Here, an experimental greenhouse is regulated via a new approach, called \emph{model-free control} (\cite{ijc13}), and their corresponding \emph{intelligent} controllers, where:
{\begin{itemize}
\item any need of a mathematical model disappears,
\item the flaws of conventional PID and On/Off techniques vanish.
\end{itemize}}
It should be emphasized that this setting (which is less than ten years old):
\begin{itemize}
\item has already been most successfully applied in a number of practical case-studies, which cover a large variety of domains (see the references in \cite{ijc13,ecc}),
\item is easy to implement (\cite{ijc13,nice}).
\end{itemize}
Besides excellent experimental results, a straightforward fault tolerant control with respect to actuators is a quite exciting byproduct. It should be emphasized here that fault accommodation for greenhouse control has
unfortunately not been very much investigated until now (see nevertheless \cite{Bontsema}).

Our paper is organized as follows. Sections \ref{free} and \ref{accom} summarize  respectively model-free control and actuator fault accommodation. Our experimental greenhouse system and its climate management problem are described in Section \ref{Apply}. Section \ref{results} displays our experimental results with our very simple intelligent controller. Comparisons with a classical Boolean controller are found in Section \ref{comp}. The efficiency of our method, is further
confirmed in Section \ref{reference} where the temperature references are modified. Section \ref{acc} deals with fault accommodation.
Some concluding remarks are provided in Section \ref{conc}.

{When compared to the two first drafts of this work, which appeared in conferences (\cite{Laflast,Laflast2}), this paper:
\begin{itemize}
\item is proposing a much simpler control synthesis than in \cite{Laflast},
\item gives a much more detailed review of model-free control than in \cite{Laflast,Laflast2},
\item reports, contrarily to \cite{Laflast,Laflast2}:
\begin{itemize}
\item the hygrometry control,
\item the time evolution of $F$ in Equation \eqref{1}.
 \end{itemize}
\end{itemize}
}

\section{Model-free control and intelligent controllers\protect\footnote{See \cite{ijc13} for more details.}} \label{free}
\subsection{The ultra-local model}\label{ulm}
For the sake of notational simplicity, let us restrict ourselves to single-input single-output (SISO) systems.\footnote{See also Section \ref{results}.} The unknown global description of the plant is replaced by the \emph{ultra-local model}:
\begin{equation}
\boxed{\dot{y} = F + \alpha u} \label{1}
\end{equation}
where:
\begin{itemize}
\item the control and output variables are respectively $u$ and $y$,
\item the derivation order of $y$ is $1$ like in most concrete situations,
\item $\alpha \in \mathbb{R}$ is chosen by the practitioner such that $\alpha u$ and
$\dot{y}$ are of the same magnitude.
\end{itemize}
The following comments might be useful:
\begin{itemize}
\item {Equation \eqref{1} is only valid during a short time lapse. It must be continuously updated,}\footnote{{The following comparison with computer graphics, which is extracted from \cite{ijc13}, might be enlightening.
Reproducing on a screen a complex plane curve is not achieved via the equations defining that curve but by approximating
it with short straight line segments. Equation \eqref{1} might be viewed as a kind of analogue of such a short segment.}}
\item $F$ is estimated via the knowledge of the control and output variables $u$ and $y$,
\item $F$ subsumes not only the unknown structure of the system, {which most of the time will be nonlinear}, but also of
any disturbance.\footnote{See also the recent comments by \cite{gao}.}
\end{itemize}
\begin{remark}
The general ultra-local model reads
$$
y^{(\nu)} = F + \alpha u
$$
where $y^{(\nu)}$ is the derivative of order $\nu \geq 1$ of $y$. When compared to Equation \eqref{1}, the only concrete case-study where such an extension was until now needed, with $\nu = 2$, has been provided by
a magnetic bearing (see \cite{compiegne}). This is explained by a very low friction (see \cite{ijc13}).
\end{remark}
\subsection{Intelligent controllers}
Close the loop with the following \emph{intelligent proportional-integral controller}, or \emph{iPI},\footnote{{The term \emph{intelligent} is borrowed from \cite{ijc13}, and from earlier papers which are cited there.}}
\begin{equation}\label{ipi}
u = - \frac{F - \dot{y}^\ast + K_P e + K_I \int e}{\alpha}
\end{equation}
where:
\begin{itemize}
\item $e = y - y^\star$ is the tracking error,
\item $K_P$, $K_I$ are the usual tuning gains.
\end{itemize}
When $K_I = 0$, we obtain \emph{intelligent proportional controller}, or \emph{iP}, which will be employed here:
\begin{equation}\label{ip}
\boxed{u = - \frac{F - \dot{y}^\ast + K_P e}{\alpha}}
\end{equation}
Combining Equations \eqref{1} and \eqref{ip} yields:
$$
\dot{e} + K_P e = 0
$$
where $F$ does not appear anymore. The tuning of $K_P$ is therefore quite straightforward. This is a major benefit when
compared to the tuning of ``classic'' PIDs (see, \emph{e.g.},
\cite{astrom,od}, and the references therein). {Note moreover that, according to Section 6.1 in \cite{ijc13}, our iP is equivalent in some sense to a classic PI controller. The integral term in the PI controllers explains why steady
state errors are avoided here with our iP.}
{\begin{remark}
Section 6 in \cite{ijc13} extends the above equivalence to classic PIDs and the
``intelligent'' controllers of \cite{ijc13}. Two important facts, which were quite mysterious in today's literature, are therefore fully clarified:
\begin{itemize}
\item the strange ubiquity of PIDs in most diverse engineering situations,
\item the difficulty of a ``good'' PID tuning for concrete industrial plants.
\end{itemize}
\end{remark}}

{\begin{remark}
Besides numerous academic comparisons in \cite{ijc13}, see, \emph{e.g.}, \cite{brest} for a thorough comparison between our intelligent controllers and PIDs for a concrete case-study, \emph{i.e.}, the position control of a shape
memory alloy active spring. All those comparisons turn out to be in favor of our intelligent controllers.
\end{remark}}


\begin{remark}
Our intelligent controllers are successfully used in an on-off way. This was also the case in \cite{sofia} for a freeway ramp metering control.
\end{remark}

\subsection{Estimation of $F$}\label{F}
Assume that $F$ in Equation \eqref{1} is ``well'' approximated by a piecewise constant function $F_{\text{est}} $. The estimation techniques below are borrowed
from \cite{sira1,sira2}.\footnote{See also the excellent recent book by \cite{sira}.}
\subsubsection{First  approach}
Rewrite then Equation \eqref{1}  in the operational domain (see, \emph{e.g.}, \cite{yosida}):
$$
sY = \frac{\Phi}{s}+\alpha U +y(0)
$$
where $\Phi$ is a constant. We get rid of the initial condition $y(0)$ by multiplying both sides on the left by $\frac{d}{ds}$:
$$
Y + s\frac{dY}{ds}=-\frac{\Phi}{s^2}+\alpha \frac{dU}{ds}
$$
Noise attenuation is achieved by multiplying both sides on the left by $s^{-2}$. It yields in the time domain the realtime estimate, thanks to the equivalence between $\frac{d}{ds}$ and the multiplication by $-t$,
\begin{equation*}\label{integral}
{\small F_{\text{est}}(t)  =-\frac{6}{\tau^3}\int_{t-\tau}^t \left\lbrack (\tau -2\sigma)y(\sigma)+\alpha\sigma(\tau -\sigma)u(\sigma) \right\rbrack d\sigma }
\end{equation*}
where $\tau > 0$ might be quite small. This integral, which is a low pass filter, may of course be replaced in practice by a classic digital filter.

\subsubsection{Second approach}\label{2e}
Close the loop with the iP \eqref{ip}. It yields:
$$
F_{\text{est}}(t) = \frac{1}{\tau}\left[\int_{t - \tau}^{t}\left(\dot{y}^{\star}-\alpha u
- K_P e \right) d\sigma \right]
$$

{\begin{remark}
It should be emphasized that the above estimation of the function $F$ in Equation \eqref{1} is quite different from model-based parameter identification. This remains valid in a control adaptive setting, where, as stated by, \emph{e.g.}, \cite{landau}, ``one needs
to know the dynamic model of the plant to be controlled.''
\end{remark}}
\begin{remark}
Implementing our intelligent controllers is easy (see \cite{ijc13,nice}).
\end{remark}

\section{Actuator's fault accommodation}\label{accom}
\begin{figure}[!t]
      \centering
      \includegraphics[scale=0.46]{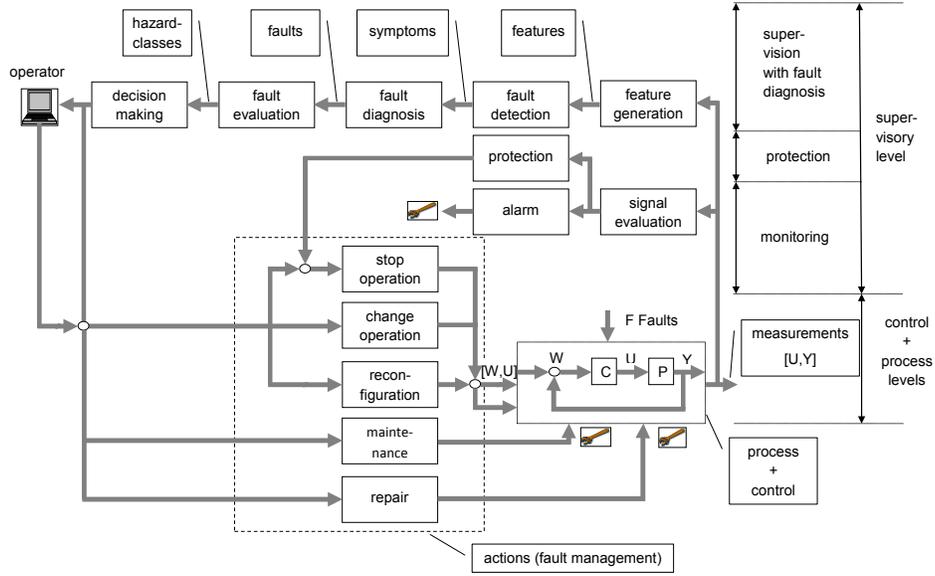}
      \caption{A supervision structure}
      \label{Sup}
 \end{figure}

As explained in Figure \ref{Sup} there are two main ways in order to deal with an actuator fault (see, \emph{e.g.}, \cite{isermann,nancy,shum}):
\begin{enumerate}
\item the first one is self-tuning, or fault accommodation. It relies on an on-line control law that preserves the main performances, while some minor parts may slightly deteriorate,
\item the second one is self-organization where faulty components are replaced.
 \end{enumerate}
 We only consider here fault accommodation. The computations below are adapted from \cite{ijc13}.

Express the actuator fault via
\begin{equation}
u_{r}  = u \left(1 - \beta\right)  \label{acc1}
\end{equation}
where:
\begin{itemize}
\item  $\beta$, $0 < \beta < 1$, is the loss of efficiency of the actuator,
\item $u_r$ is the true control variable.
\end{itemize}
The two following cases are not considered:
\begin{itemize}
\item  $\beta = 0$ means that there is no fault,
\item $\beta = 1$ implies that the control does not act anymore.
\end{itemize}
Then Equation \eqref{1} becomes
$$
\dot{y} = \bar{F} + \alpha u
$$
where
$$
\bar{F} = F - \alpha \beta u
$$
The fault accommodation is then achieved by estimating $\bar{F}$ as in Section \ref{F}.

\begin{remark}
It is obvious that $\beta$ does not need to be:
\begin{itemize}
\item a constant and may be time-varying,
\item known in order to carry on the above computations.
 \end{itemize}
\end{remark}

\begin{remark}
For model-based diagnosis, estimation techniques stemming from \cite{sira1,sira2} have already lead to quite important advances. See, \emph{e.g.}, \cite{ijc04,nl,sarre,mines1,mines2}.
\end{remark}

\section{Greenhouse climate management} \label{Apply}
\begin{figure}[!tbp]
      \centering
      \includegraphics[scale=0.8]{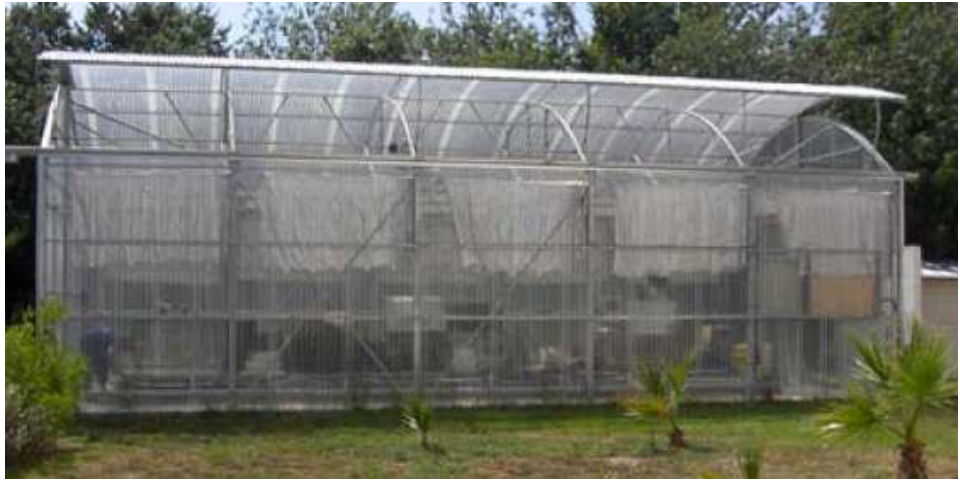}
      \caption{Our experimental greenhouse system}
      \label{Greenhouse}
   \end{figure}

Figure \ref{Greenhouse} shows our experimental plastic greenhouse which is manufactured by the French company \emph{Richel}. Its area is equal to $80$ m$^{2}$. It is the property of the \emph{Laboratoire des Sciences de l'Information
et des Syst\`{e}mes}  (\emph{LSIS}), to which the first three authors belong. This laboratory is located at the \emph{Universit\'e de Toulon} in the south of France. 
Our experimental greenhouse is controlled by a microcomputer and interfaced with the FieldPoint FP-2000 network module
developed by the American company \emph{National Instruments Corporation}. The FP-2000 network module is associated with two analog input modules (FP-AI-110, FP-AI-111), for the acquisition, and two
relay output modules (FP-RLY-420), for the control. The acquisition and control system is developed with the \emph{LabView} language. The sampling period is equal to $1$ minute.
The inside air temperature and the humidity are controlled.

\subsection{Description of the system}
\begin{figure}[!tbp]
      \centering
      \includegraphics[scale=1]{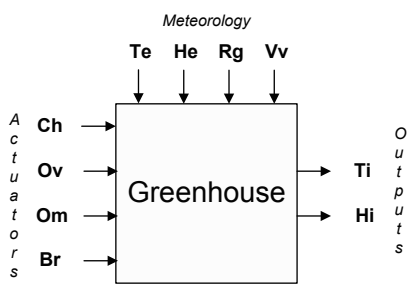}
      \caption{System variables}
      \label{Variables}
   \end{figure}

The greenhouse is a multi-input and multi-output (MIMO) system which is equipped with several sensors and actuators (Figure \ref{Variables}). \\

\begin{figure}[!tbp]
      \centering
      \includegraphics[scale=0.6]{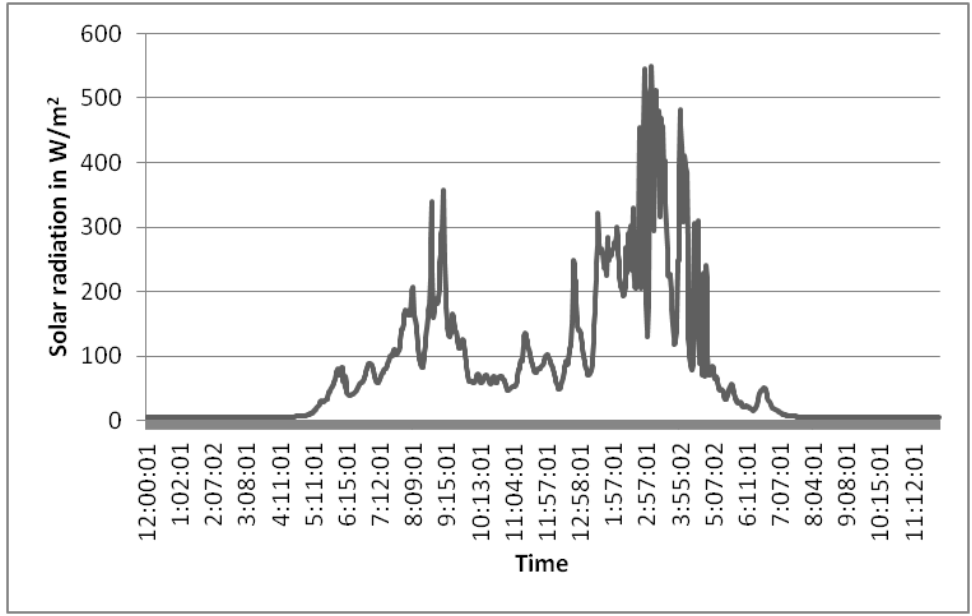}
      \caption{Solar radiation during the $24^{th}$ September, 2014}
      \label{Variable2}
   \end{figure}

\begin{figure}[!tbp]
      \centering
      \includegraphics[scale=0.6]{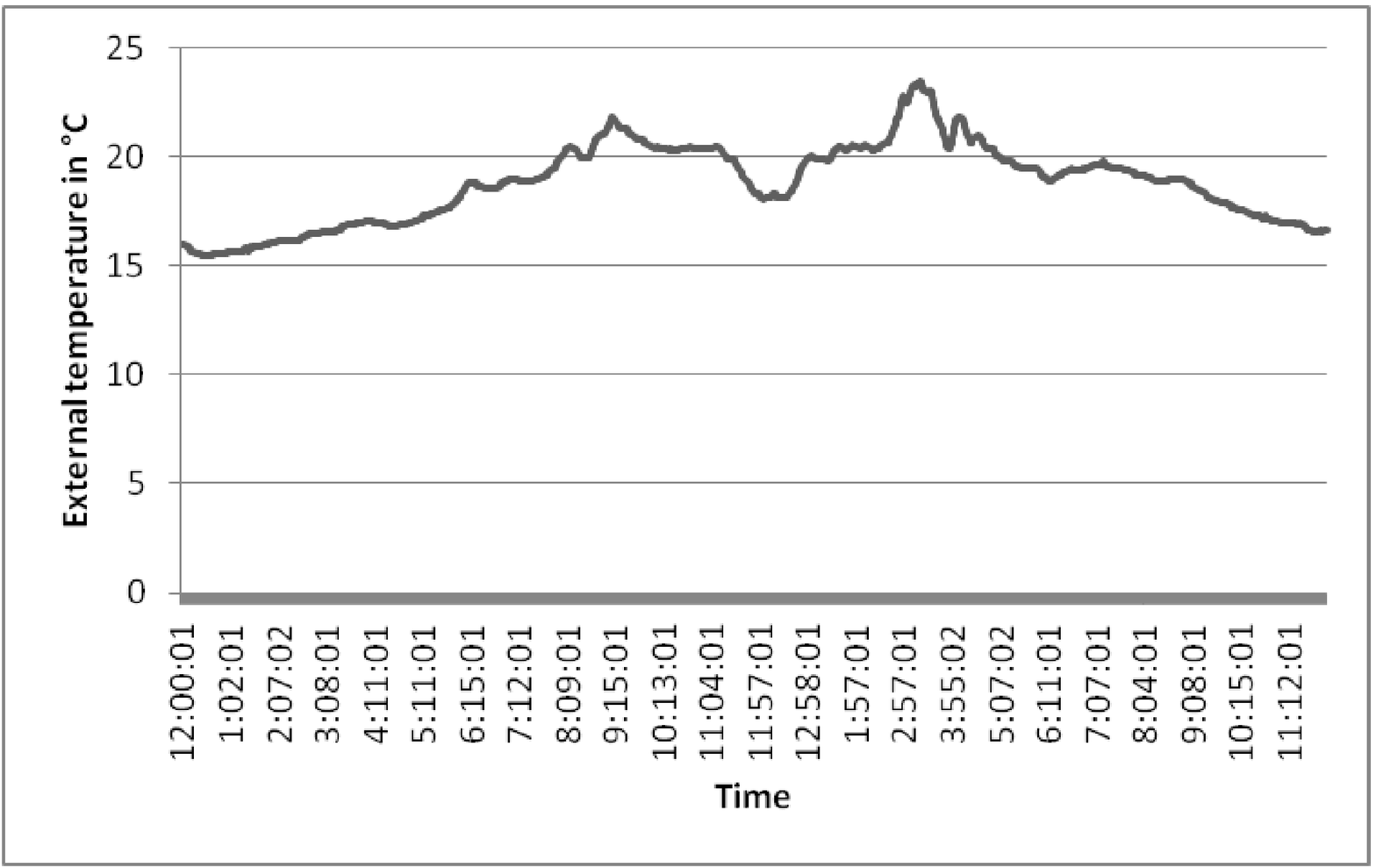}
      \caption{External temperature during the $24^{th}$ September, 2014}
      \label{Variable3}
   \end{figure}

\begin{figure}[!tbp]
      \centering
      \includegraphics[scale=0.6]{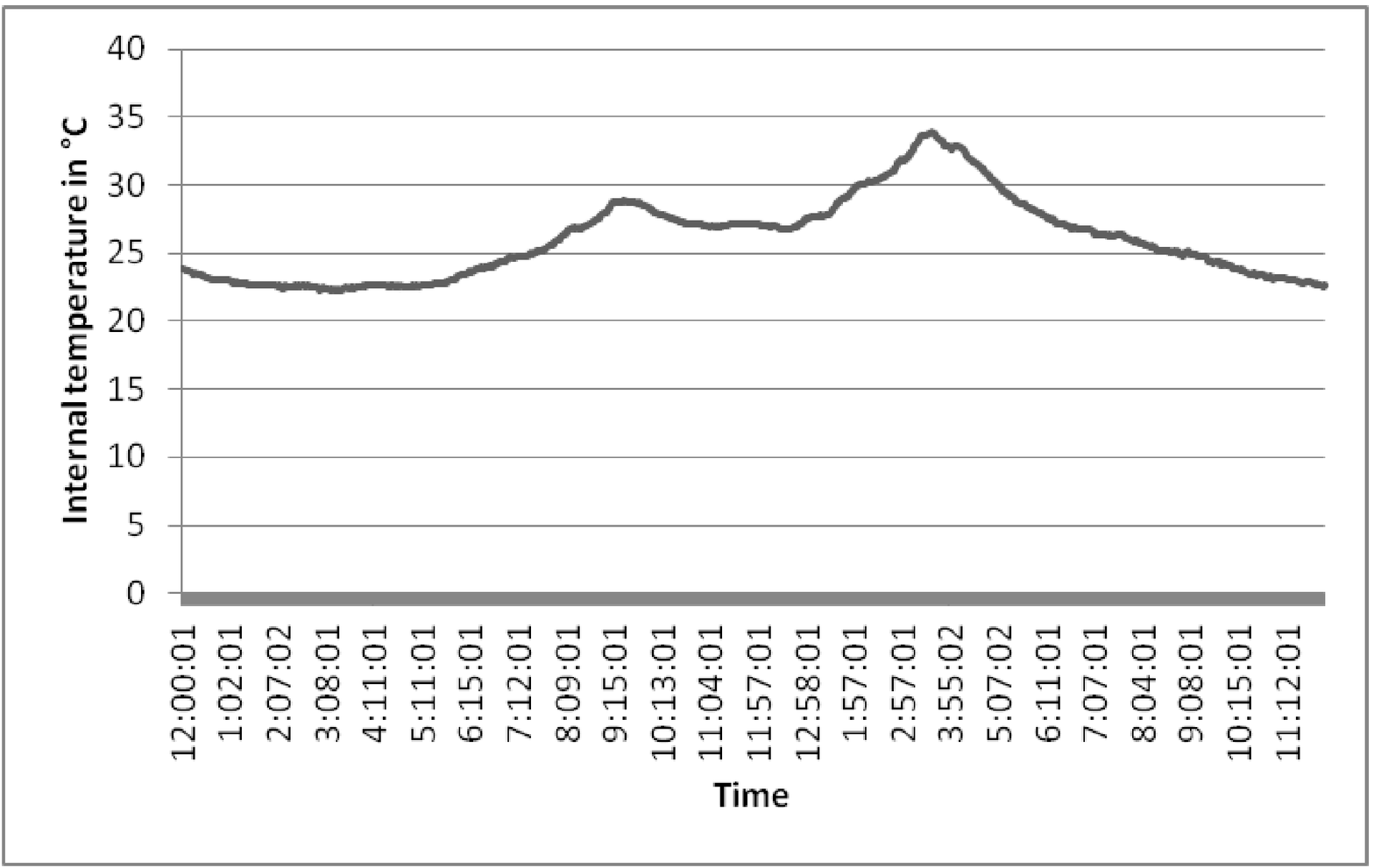}
      \caption{Internal temperature during the $24^{th}$ September, 2014}
      \label{Variable4}
   \end{figure}

There are:
\begin{itemize}
\item four actuators:
\begin{enumerate}
 \item Heating (thermal power 58 kw): Ch ($Boolean$),
 \item Opening (50 \% max): Ov (\%),
 \item Shade: Om (\%),
 \item Fog system: Br ($Boolean$).
\end{enumerate}

\item four meteorological disturbance sensors:
\begin{enumerate}
 \item External temperature: Te ($^{o}C$),
 \item External hygrometry: He (\%),
 \item Solar Radiation: Rg ($W/m^{2}$),
 \item Wind speed: Vv ($km/h$).
 \end{enumerate}

\item two internal climate sensors:
\begin{enumerate}
 \item Internal temperature: Ti ($^{o}C$),
  \item Internal hygrometry: Hi (\%).\\
  \end{enumerate}
\end{itemize}
\noindent This system is nonstationary and strongly disturbed. Figures \ref{Variable2} and \ref{Variable3} show, for instance, quite high solar radiation and external temperature during the $24^{th}$ September 2014.
These meteorological conditions have a significant effect on the inside greenhouse climate which are clear on Figure \ref{Variable4}.

\subsection{Climate management problem} \label{Refer}
The management of the greenhouse climate aims to maintain simultaneously a set of climatic factors such as the temperature, the hygrometry, and the rate of $CO_{2}$\footnote{This last rate is not available on our greenhouse.} close to their respective references.
In our greenhouse, the temperature and the hygrometry managements are treated together, because these two
quantities are strongly correlated:
\begin{itemize}
\item the heating has a dehumidifier effect,
\item the opening system has a cooling and dehumidifier effect,
\item the fog system has a cooling effect.
\end{itemize}
Controlling the temperature and the hygrometry is therefore of utmost importance. In order to choose the suitable output references, two main strategies exist.

\subsubsection{The classic strategy} \label{clastra}
Growers refer to their knowledge to fix the hygrometry and temperature references.

\paragraph{Hygrometry reference}
There is no real recommendations by species. It appears nevertheless that:
\begin{itemize}
\item for the multiplication phase, the hygrometry must be greater than 80 \%,
\item for the growth phase, the reference is comprised between 60 and 80 \%,
\item for the tomato, the reference is rather comprised between 50 and 70 \%.\\
\end{itemize}
Let us mention some other advices. Avoid:
\begin{itemize}
\item condensations,
\item a humidity level close to saturation (100 \%),
\item a humidity level below 40 \% for seedlings,
\item absolutely a hygrometry below 20 \%.
\end{itemize}

\paragraph{Temperature reference}
Table \ref{Reference} displays references among suppliers, which are based on the species.\footnote{Temperatures are expressed in Celsius degrees.}
Observe that the difficulties for tuning an efficient controller may be attributed to the following causes:
\begin{itemize}
\item various references:
\begin{itemize}
\item in a day,
\item according to the species.
\end{itemize}
\item system parameter variations according to the plant growth.\\
\end{itemize}

\begin{table}[!tbp]
\caption{Temperature reference (see \citep{c15})}
\label{Reference}
\begin{center}
\begin{tabular}{|c||c||c||l|}
\hline
Species & Night & Day & Remarks\\
 & reference & reference & \\
\hline
Aubergine & 21$^{o}C$  & 22$^{o}C$ & During 4 weeks\\
& & & after the plant.\\
 & 19$^{o}C$ & 21$^{o}C$ & To the end\\
\hline
Cucumber & 21$^{o}C$  & 23$^{o}C$ & During 4 weeks\\
& & & After the plant.\\
 & 20$^{o}C$ & 22$^{o}C$ & During the next\\
& & & 6 weeks.\\
 & 19$^{o}C$ & 21$^{o}C$ & To the end.\\
 \hline
Lettuce & 10$^{o}C$  & 10$^{o}C$ & During 2 weeks\\
& & & After the plant.\\
 & 6$^{o}C$ & 12$^{o}C$ & To the end.\\
 \hline
Pepper & 20$^{o}C$  & 23$^{o}C$ & During 3 weeks\\
& & & after the plant.\\
 & 18$^{o}C$ & 22$^{o}C$ & To the end.\\
 \hline
Tomato & 20$^{o}C$  & 20$^{o}C$ & During 1 week\\
& & & after the plant.\\
 & 18.5$^{o}C$ & 19.5$^{o}C$ & During the next\\
& & & 5 weeks.\\
 & 17.5$^{o}C$ & 18.5$^{o}C$ & To the end.\\
 \hline
Azalea & 18/21$^{o}C$  & $>$18$^{o}C$ & \\
\hline
Chrysanthemum & 17$^{o}C$  & 18$^{o}C$ & \\
\hline
Gerbera & 13/15$^{o}C$  &  & \\
\hline
Antirrhinum & 10/11$^{o}C$  &  & \\
\hline
Carnation & 12/13$^{o}C$  & 18$^{o}C$ & \\
\hline
Rosebush & 17$^{o}C$  & 21$^{o}C$ & \\
\hline
\end{tabular}
\end{center}
\end{table}

\subsubsection{The innovative strategy} \label{innstra}
\cite{tcham} developped a decision-making system, called SERRISTE. It generates daily climate reference for greenhouse grown tomatoes. This system, which uses the knowledge of advisers or expert growers to manage the greenhouse climate, can be encapsulated and exploited in a reference determination software. This tool provides daily references to growers taking into account various objectives such as the phytosanitary prevention, the energetic cost, the growth of the crop, ... . The system uses data such as seasons, crop stages, the daily period (divided into three subperiods), the characteristics of the greenhouse system (location, heating system, ...) and dynamic informations (past climate, crop state, ...).
Sections \ref{clastra} and \ref{innstra} show the reference changes according to the time of day or the plant growth. This is another justification for our model-free control.

\section{Intelligent P control of the experimental greenhouse} \label{results}
\begin{figure}[!tbp]
      \centering
      \includegraphics[scale=0.5]{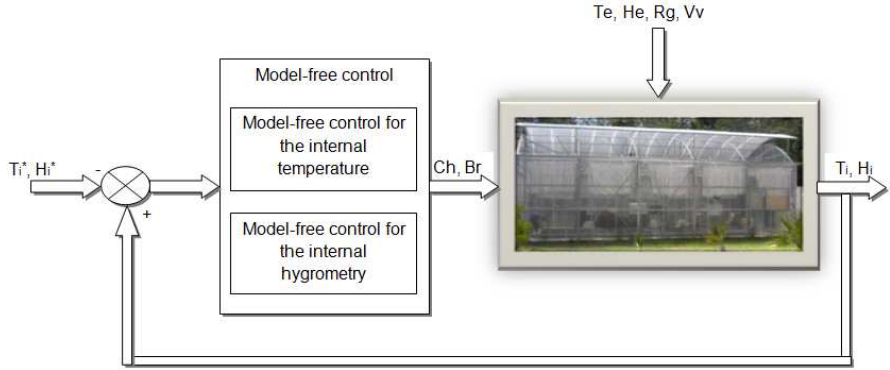}
      \caption{Block diagram of the experimental setup}
      \label{BF}
   \end{figure}

An iP \eqref{ip} is implemented for the regulation of the temperature and the hygrometry, which turn out to be naturally decoupled in our model-free setting (Figure \ref{BF}).\footnote{Our restriction in Section \ref{free} to detail only SISO systems is therefore
fully justified. See also \cite{menhour} for the behavior of a vehicle.}

We are estimating $F$ via the technique sketched in Section \ref{2e}.

\subsection{Estimation of $F$}
The estimation $F^{\text{temp}}_{\text{est}}$ is given by
\begin{equation}\label{eq12}
F^{\text{temp}}_{\text{est}} = \frac{1}{\delta} \int^{T}_{T-\delta} \left(- \alpha Ch + \dot{Ti}^{*} - K_{P} e_{Ti}\right) d\tau
\end{equation}
where:
\begin{itemize}
\item $e_{Ti} = Ti - Ti^{*}$ is the temperature tracking error,
\item $\dot{Ti}^{*}$ is the reference derivative of $Ti$ (when internal temperature reference is constant then $\dot{Ti}^{*}$ is equal to 0).
\end{itemize}

and $F^{\text{hygro}}_{\text{est}}$ by
\begin{equation}\label{eq13}
\begin{array}{c}
F^{\text{hygro}}_{\text{est}} = \frac{1}{\delta} \int^{T}_{T-\delta} \left(- \alpha Br + \dot{Hi}^{*} - K_{P} e_{Hi}\right) d\tau
\end{array}
\end{equation}

where:
\begin{itemize}
\item $e_{Hi} = Hi - Hi^{*}$ is the temperature tracking error,
\item $\dot{Hi}^{*}$ is the reference derivative of $Hi$ (when internal hygrometry reference is constant then $\dot{Hi}^{*}$ is equal to 0).
\end{itemize}

\subsection{Setting values and results}
The controllers $Ch$ and $Br$  are deduced from Equations \eqref{1}, \eqref{ip} and \eqref{eq12}.
They are \emph{Pulse Width Modulation} (\emph{PWM}) controllers.
The rules given in Section \ref{ulm} yield Table \ref{Setting_values}, which displays the same values for the two controllers. 

\begin{table}[!tbp]
\caption{Setting values}
\label{Setting_values}
\begin{center}
\begin{tabular}{|c||c|}
\hline
Variable & Value\\
\hline
$\delta$ & 6 minutes\\
\hline
$\alpha$ & 1\\
\hline
$K_{P}$ & 2\\
\hline
\end{tabular}
\end{center}
\end{table}

The reference output is $18^{o}C$ for the temperature with a tolerance equal to $0.5^{o}C$ and 60 \% for the hygrometry. {The temperature sensors PT100 sensors, of class A, with an accuracy of $\pm$ $0.3^{o}C$. A tolerance of $0.5^{o}C$ would be realistic since, for many species, the difference between night and day reference is equal to $1^{o}C$, as shown in Table \ref{Reference}. We want to differentiate night and day. Sensors with an accuracy of $\pm$ $0.3^{o}C$ permit to take into account a tolerance equal to $0.5^{o}C$.} Simulations last 12 hours, from 8:00 p.m. until 8:00 a.m.
We choose the night in order to compare the obtained results with Boolean control (see Section \ref{comp}) in similar weather conditions.

\begin{figure}[!tbp]
      \centering
      \includegraphics[scale=0.6]{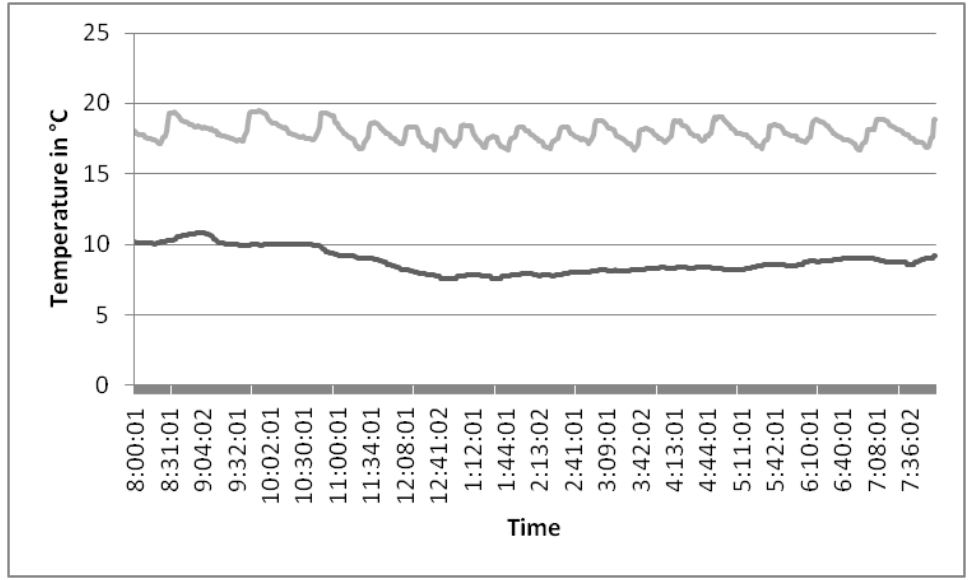}
      \caption{Temperature with model-free control (Te: black line - Ti: grey line)}
      \label{Temperature}
   \end{figure}

\begin{figure}[!tbp]
      \centering
      \includegraphics[scale=0.6]{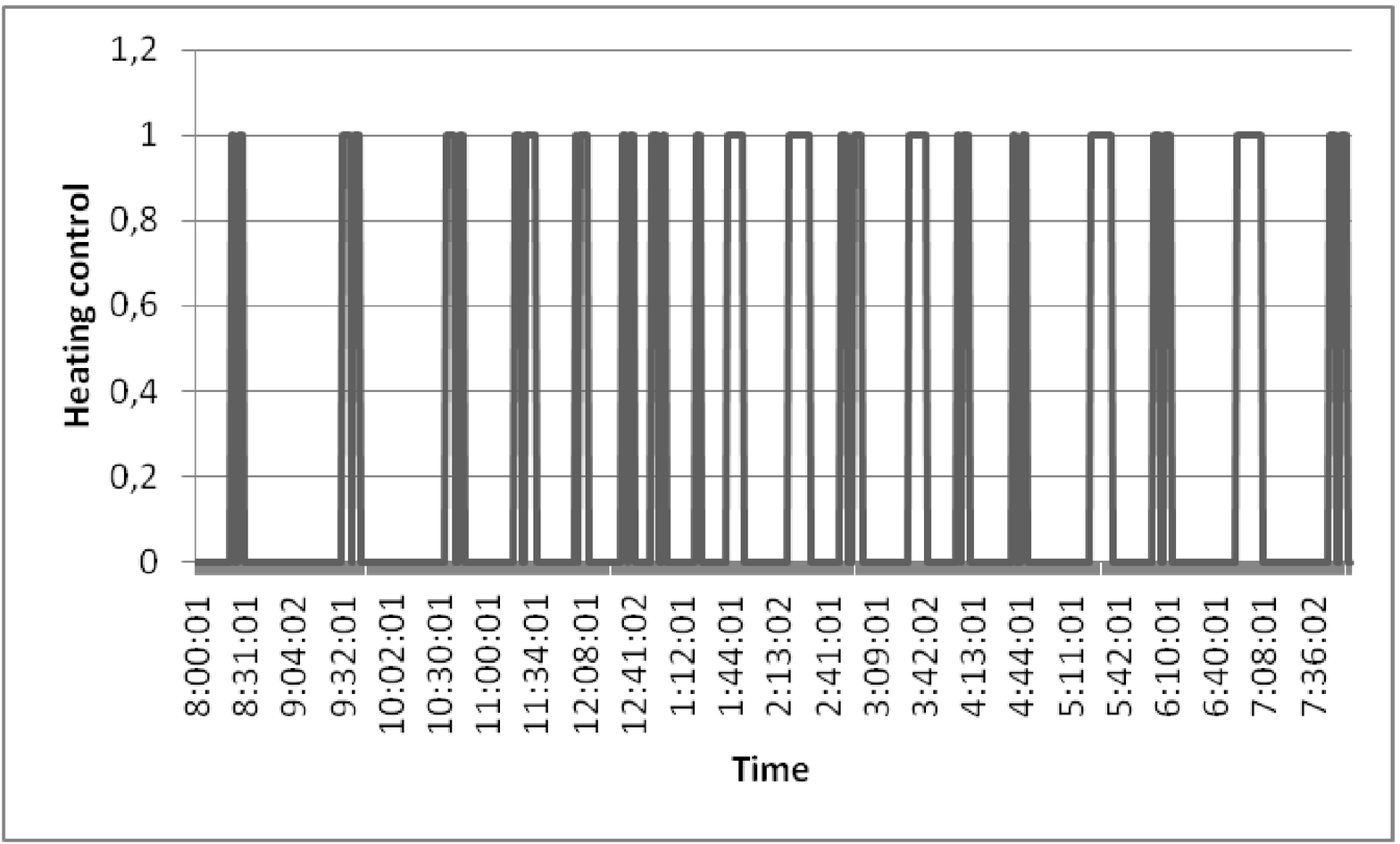}
      \caption{Heating control with model-free control}
      \label{Heating}
   \end{figure}

\begin{figure}[!tbp]
      \centering
      \includegraphics[scale=0.6]{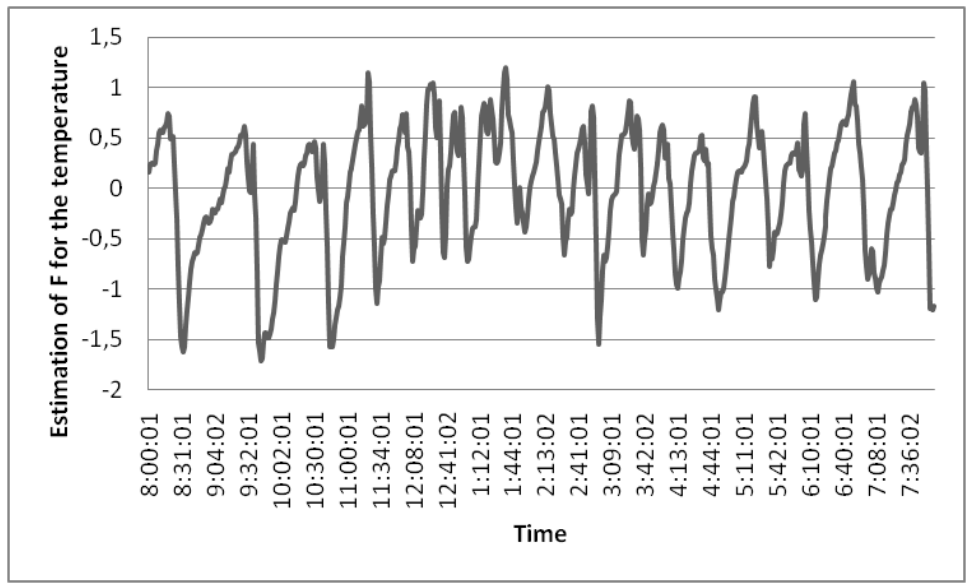}
      \caption{Evolution of $F^{temp}_{est}$}
      \label{Fest}
   \end{figure}

Figure \ref{Temperature} shows the internal/external temperature evolution during the night of 20-21 February 2014. Figure \ref{Heating} shows the heating control sequences.
Observe that the heating control allows at the internal temperature $T_{i}$ to be close to its reference output. Figure \ref{Fest} shows the evolution of {$F^{\text temp}_{\text est}$} during this night.

\begin{figure}[!tbp]
      \centering
      \includegraphics[scale=0.6]{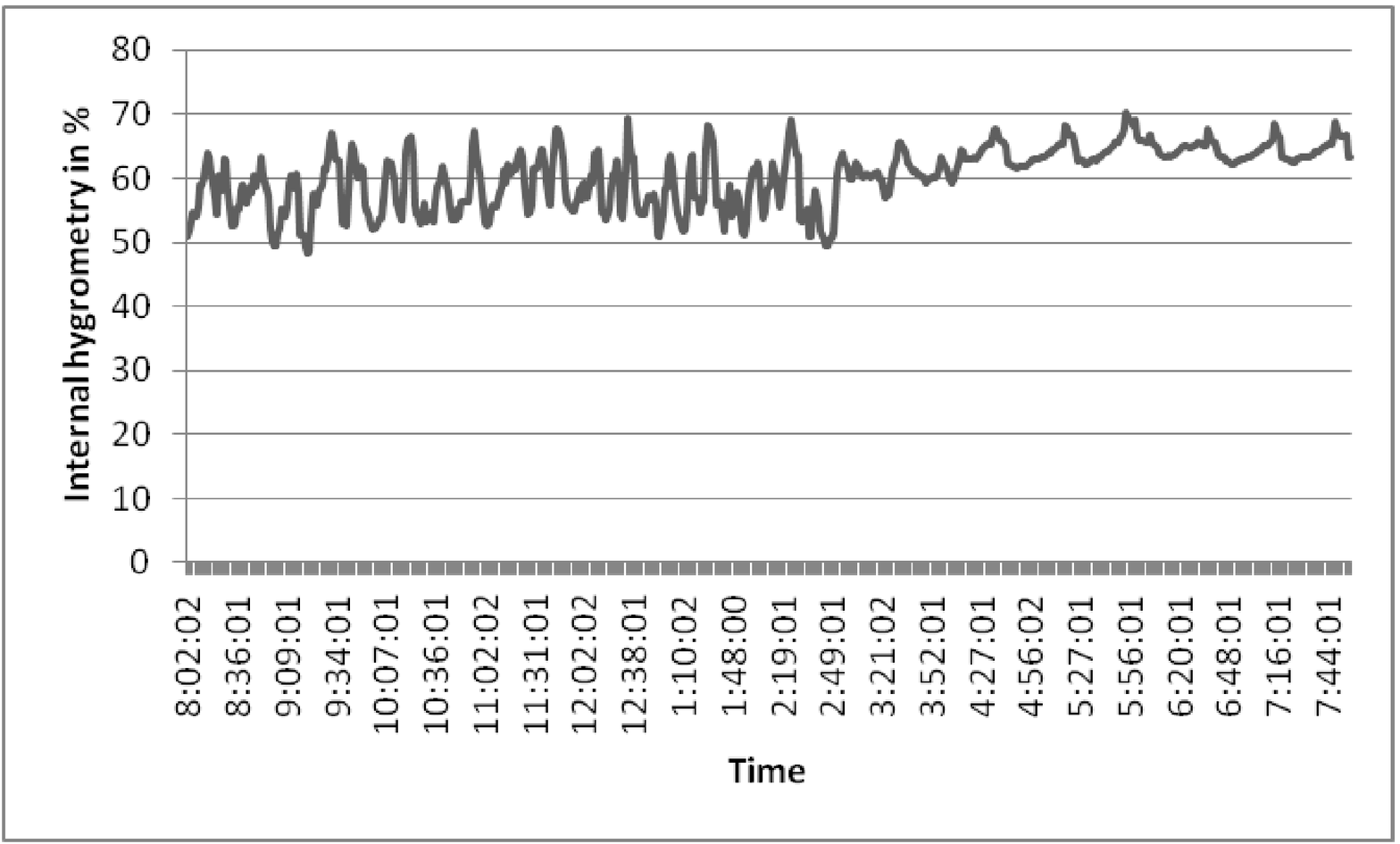}
      \caption{Internal hygrometry with model-free control}
      \label{Hygrometry}
   \end{figure}

\begin{figure}[!tbp]
      \centering
      \includegraphics[scale=0.6]{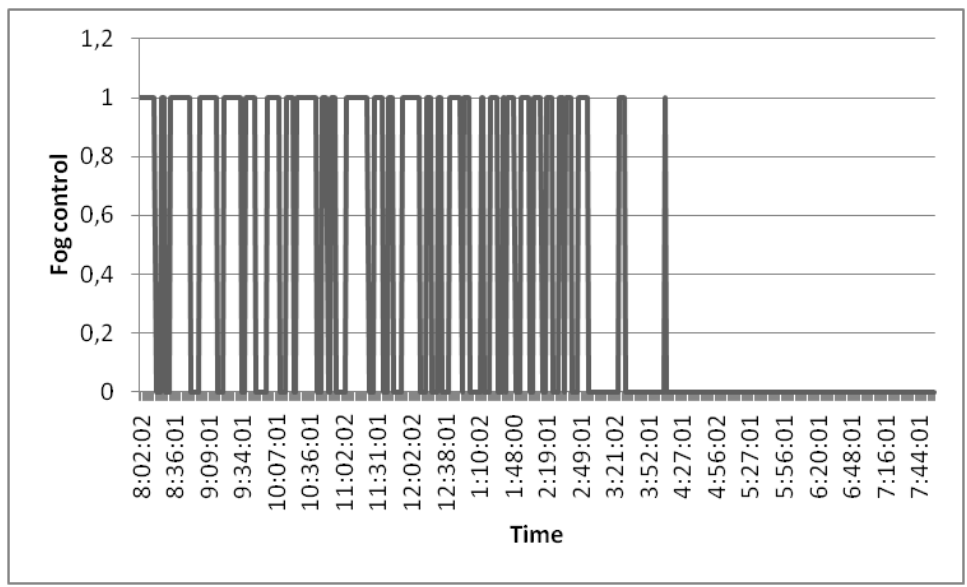}
      \caption{Fog control with model-free control}
      \label{Mois}
   \end{figure}

Figure \ref{Hygrometry} shows the internal hygrometry evolution during the night of 20-21 February 2014. Figure \ref{Mois} shows the sequences for the fog control. We can observe that, at 4:00 a.m., the internal hygrometry $H_{i}$ is also above the reference output: it started to rain. So, the fog system $Br$ stops. Otherwise, the internal hygrometry $H_{i}$ is close to this reference output.

Table \ref{Evaluation} shows the mean and the variance of the error between $T_{i}$ and the output reference of $T_{i}$ and between $H_{i}$ and the reference output of $H_{i}$.

\begin{table}[!tbp]
\caption{Results evaluation for the model-free control}
\label{Evaluation}
\begin{center}
\begin{tabular}{|c||c||c|}
\hline
Output error & mean & variance\\
\hline
$e_{T_{i}}$ & $-0.1^{o}C$ & $0.4^{o}C$\\
\hline
$e_{H_{i}}$ & 0.4 \% & 21.8 \%\\
\hline
\end{tabular}
\end{center}
\end{table}

\section{Comparison between iP and classic Boolean control} \label{comp}
A classic Boolean control law with thresholds is employed for the comparisons. This type of technique is quite often utilized in agriculture.
Experiments have been carried on during two different nights, \emph{i.e.},  20 -21 and 21-22 February 2014, respectively for the model-free and boolean settings.
The temperature reference output is $18^{o}C$ with a tolerance equal to $0.5^{o}C$, as in Section \ref{results}. For the hygrometry, a dehumidification reference should be selected. The fog control is periodic (3 minutes on and 27 minutes off) whatever the internal hygrometry. This Boolean control of the humidity is based on the grower rules. The dehumidification reference allows to set the desired maximum hygrometry inside the greenhouse. In this test, we choose 60 \%.

\begin{figure}[!tbp]
      \centering
      \includegraphics[scale=0.6]{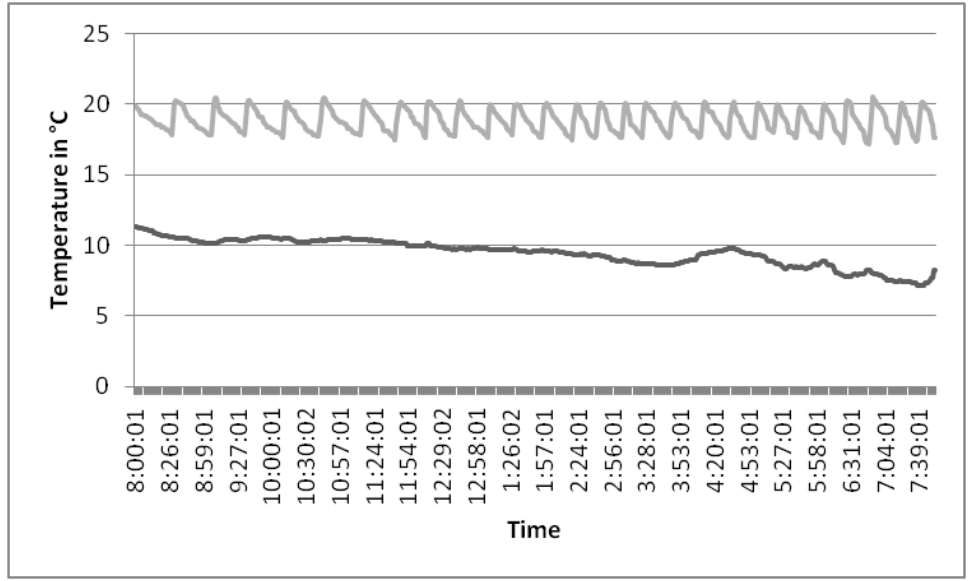}
      \caption{Temperature with a Boolean controller (Te: Black line - Ti: Grey line)}
      \label{Temperature1}
   \end{figure}

\begin{figure}[!tbp]
      \centering
      \includegraphics[scale=0.6]{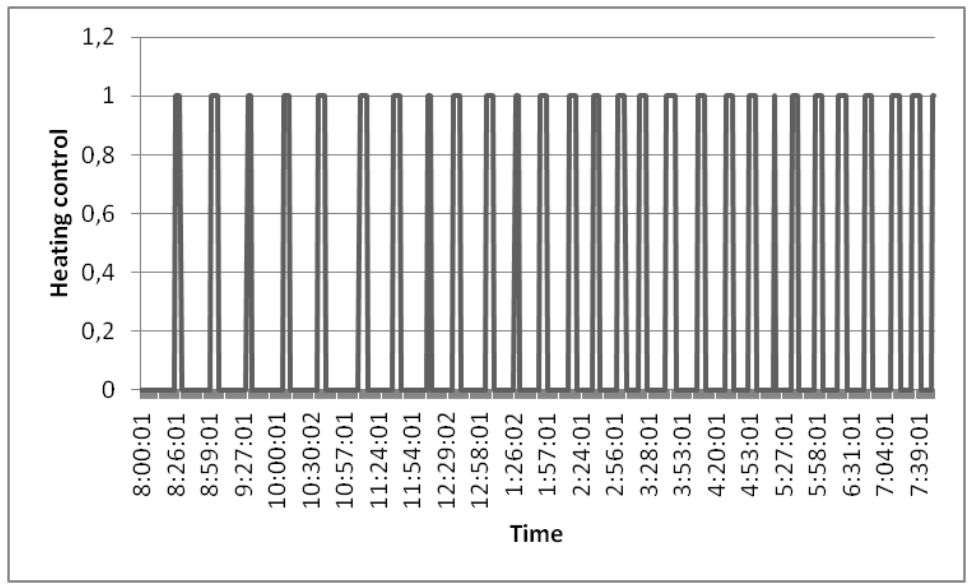}
      \caption{Heating control with a Boolean controller}
      \label{Heating1}
   \end{figure}

Figure \ref{Temperature1} and \ref{Heating1} show respectively results for the internal temperature and for the heating control during the night of 21-22 February 2014.

\begin{figure}[!tbp]
      \centering
      \includegraphics[scale=0.6]{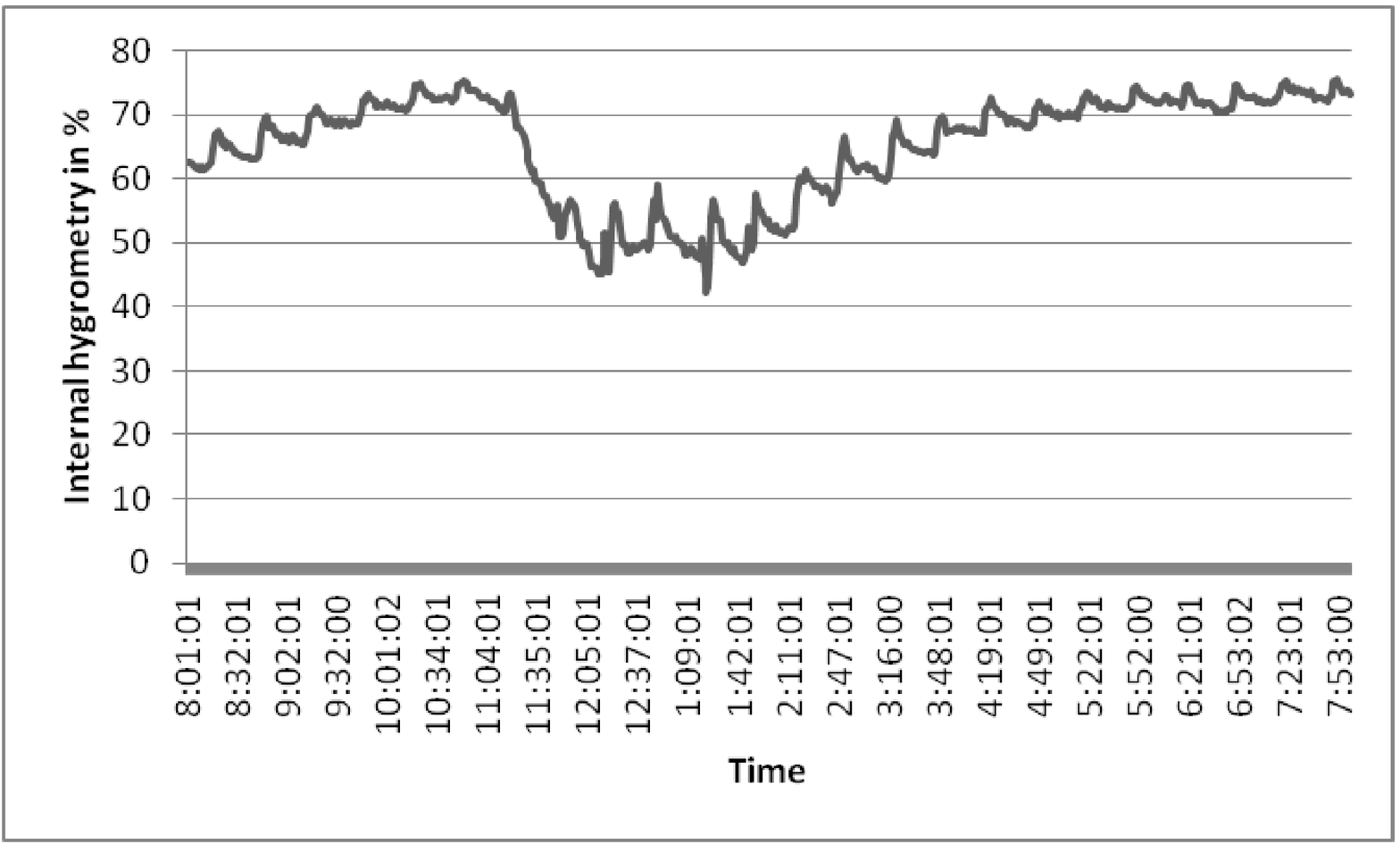}
      \caption{Internal hygrometry with a Boolean controller}
      \label{Hygrometry2}
   \end{figure}

\begin{figure}[!tbp]
      \centering
      \includegraphics[scale=0.6]{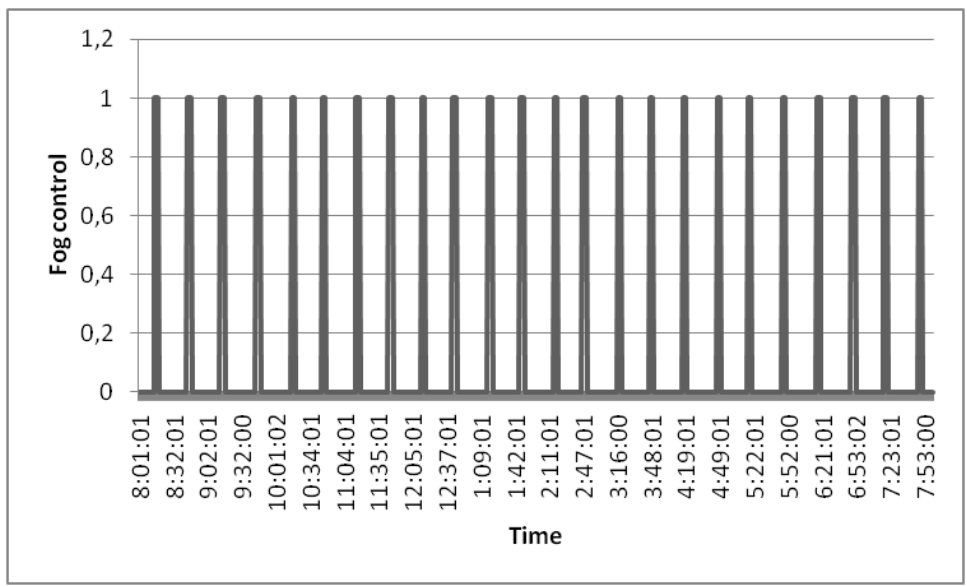}
      \caption{Fog control with a Boolean controller}
      \label{Mois2}
   \end{figure}

Figure \ref{Hygrometry2} shows the internal hygrometry evolution during the night of 21-22 February 2014. Figure \ref{Mois2} shows the sequences for the fog control.

 Table \ref{Evaluation2} shows the mean and the variance of the error between $T_{i}$ and the output reference of $T_{i}$ for this night.
\begin{table}[!tbp]
\caption{Results evaluation with a classic Boolean control}
\label{Evaluation2}
\begin{center}
\begin{tabular}{|c||c||c|}
\hline
Output error & mean & variance\\
\hline
$e_{T_{i}}$ & $0.8^{o}C$ & $0.7^{o}C$\\
$e_{H_{i}}$ & 5.0 \% & 71.7 \%\\
\hline
\hline
\end{tabular}
\end{center}
\end{table}

\begin{table}[!tbp]
\caption{Comparisons of the energy}
\label{EnergyCompare}
\begin{center}
\begin{tabular}{|c||c||c|}
\hline
Actuator & Model-free control & Classical Boolean control\\
\hline
$Heat$ & 143 $min$ & 145 $min$\\
\hline
\end{tabular}
\end{center}
\end{table}

Tables  \ref{Evaluation} and  \ref{Evaluation2} demonstrate that our model-free control strategy behaves better than its Boolean counterpart.
Let us emphasize two more points:
\begin{itemize}
\item as already explained in Section \ref{Apply}, one of the goals of climate control is to consume as little energy as possible.
Table \ref{EnergyCompare} shows that the heating is on only during 20 \% of the time with the model-free setting.
The model-free controller is therefore much cheaper,
\item for a given operating time, the model-free control ensures a better tracking of the reference signal.
\end{itemize}

\section{Reference change} \label{reference}
\begin{figure}[!tbp]
      \centering
      \includegraphics[scale=0.6]{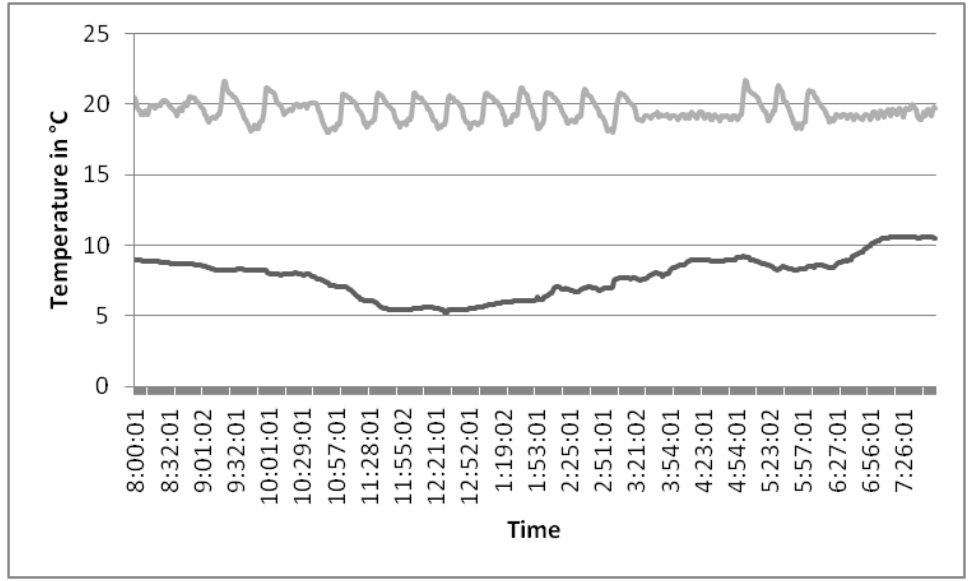}
      \caption{Temperature with model-free control (Te: Black line - Ti: Grey line)}
      \label{Temperature2}
   \end{figure}

\begin{figure}[!tbp]
      \centering
      \includegraphics[scale=0.6]{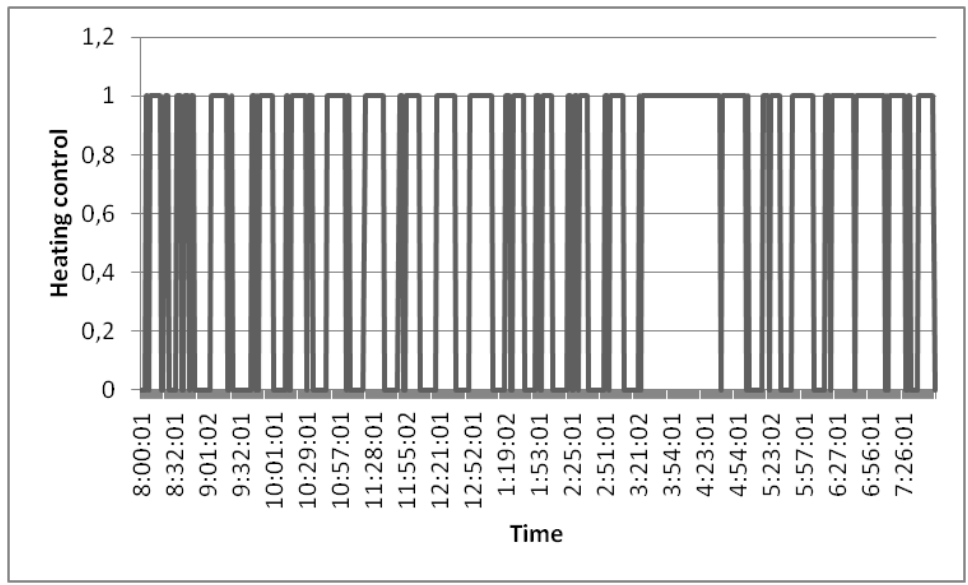}
      \caption{Heating control with model-free control}
      \label{Heating2}
   \end{figure}

Figure \ref{Temperature2} shows results for the internal temperature with a reference change (without any modification of the parameter values of the iP controller). We regulate the greenhouse with the temperature reference output equal to $20^{o}C$ during the night of 11-12 February 2014. Figure \ref{Heating2} represents the heating control.

\begin{figure}[!tbp]
      \centering
      \includegraphics[scale=0.6]{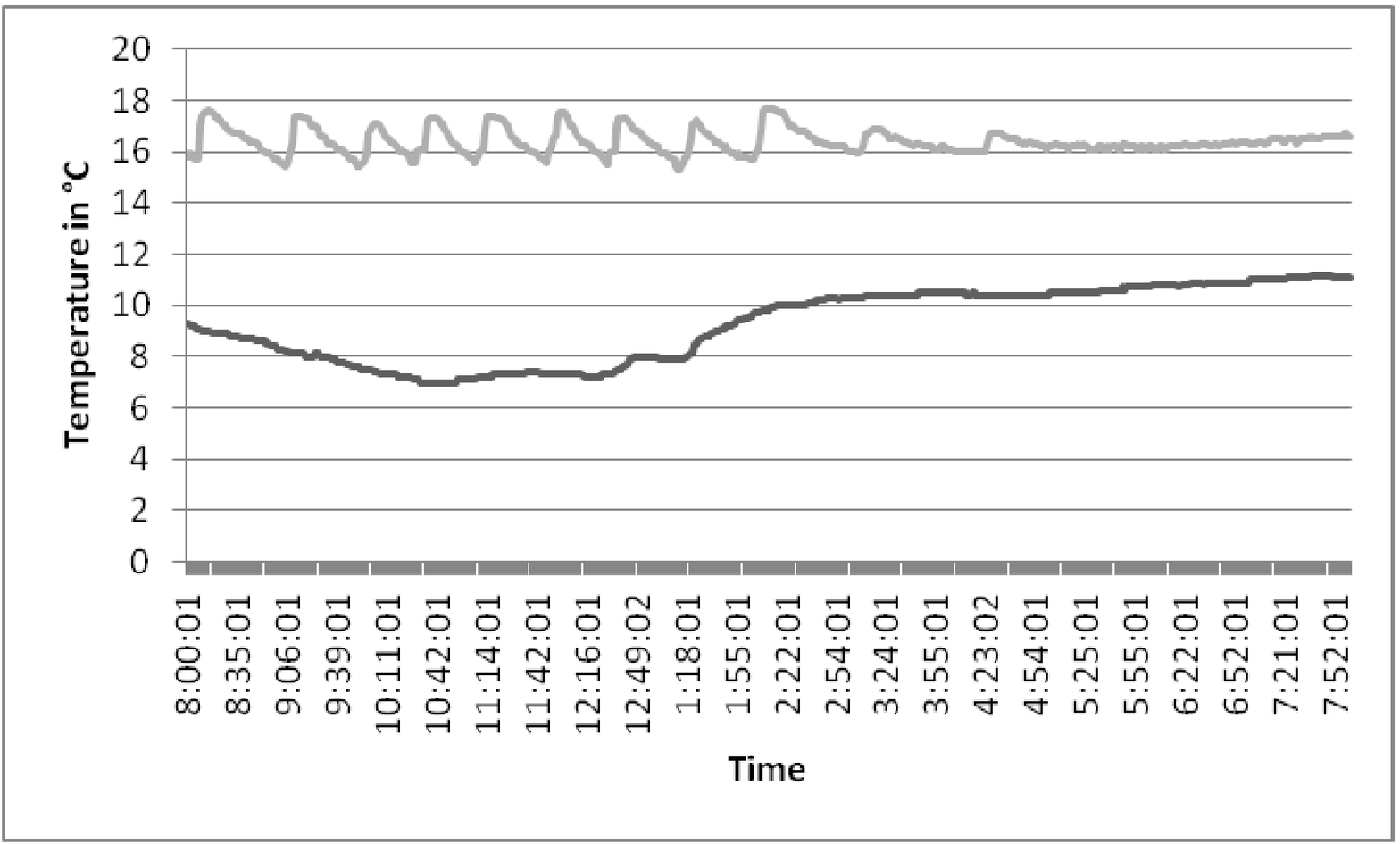}
      \caption{Temperature with model-free control (Te: Black line - Ti: Grey line)}
      \label{Temperature3}
   \end{figure}

\begin{figure}[!tbp]
      \centering
      \includegraphics[scale=0.6]{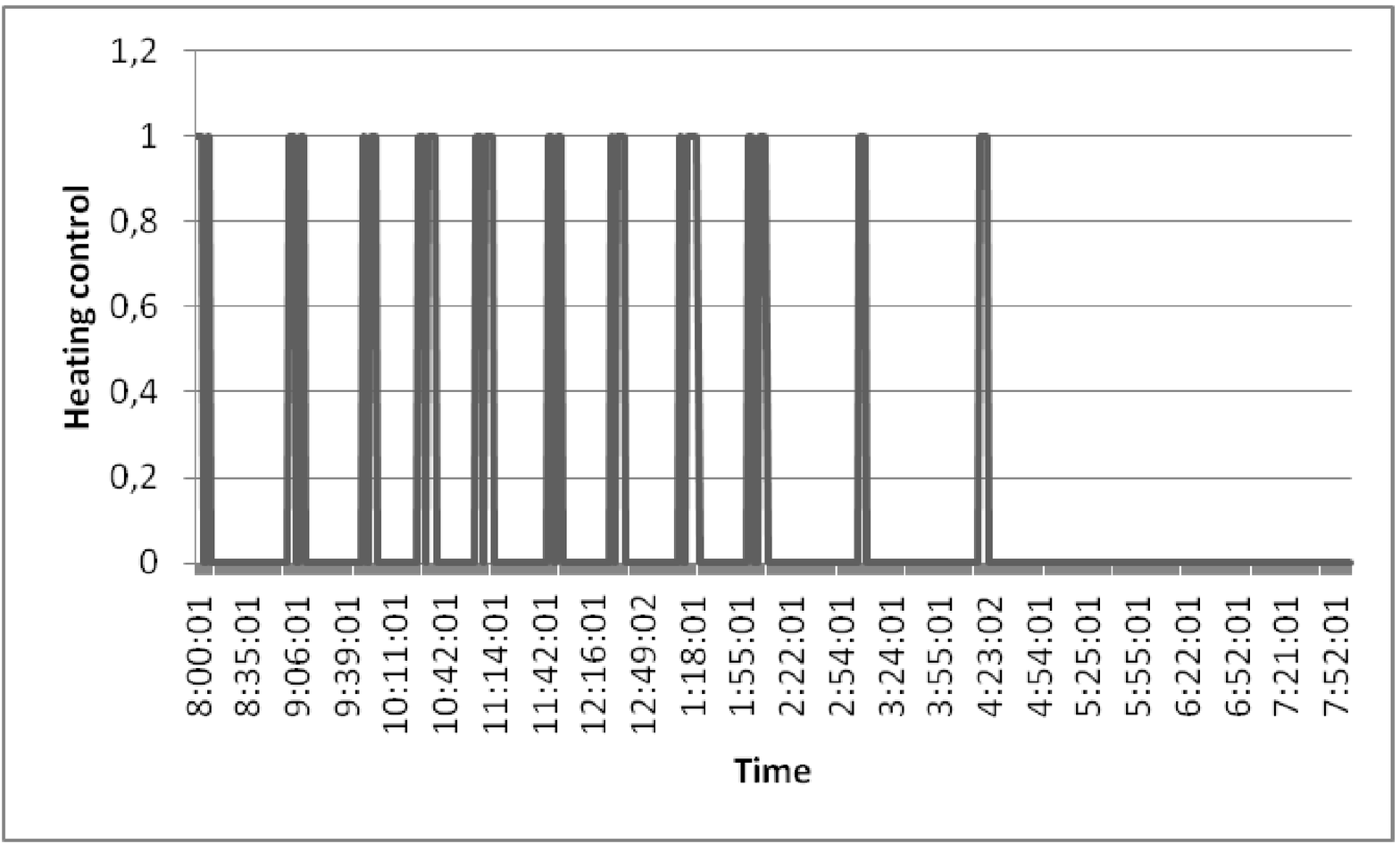}
      \caption{Heating control with model-free control}
      \label{Heating3}
   \end{figure}

Results for the internal temperature with an other reference change are displayed on Figure \ref{Temperature3}. We regulate the greenhouse with the temperature reference output equal to $16^{o}C$ during the night of 17-18 February 2014. Figure \ref{Heating3} represents the heating control.

We can observe that model-free control results are always \emph{good} since the internal temperature follow to the reference output (see Table \ref{Evaluation3}). As sketched in Section \ref{Refer} and presented in Table \ref{Reference}, this is a
most significant advance.

\begin{table}[!tbp]
\caption{Results evaluation for the model-free control}
\label{Evaluation3}
\begin{center}
\begin{tabular}{|c||c||c|}
\hline
Output error & mean & variance\\
\hline
$e_{T_{i}}$ for $Ti^{*} = 20^{o}C$ & $-0.4^{o}C$ & $0.6^{o}C$\\
\hline
$e_{T_{i}}$ for $Ti^{*} = 16^{o}C$ & $0.4^{o}C$ & $0.2^{o}C$\\
\hline
\end{tabular}
\end{center}
\end{table}

\section{Fault accommodation} \label{acc}
\begin{figure}[!tbp]
      \centering
      \includegraphics[scale=0.6]{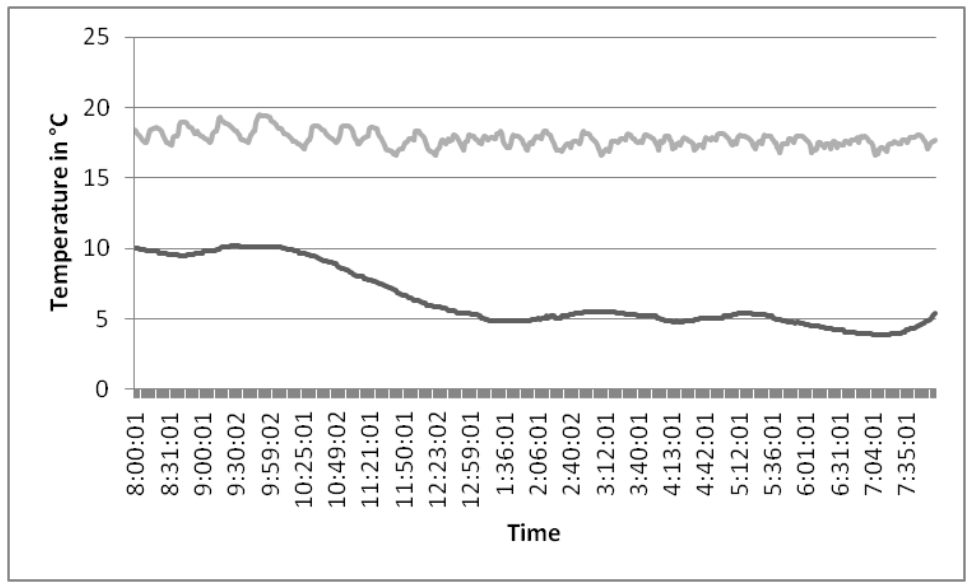}
      \caption{Temperature with model-free control (Te: Black line - Ti: Grey line)}
      \label{Temperature4}
   \end{figure}

\begin{figure}[!tbp]
      \centering
      \includegraphics[scale=0.6]{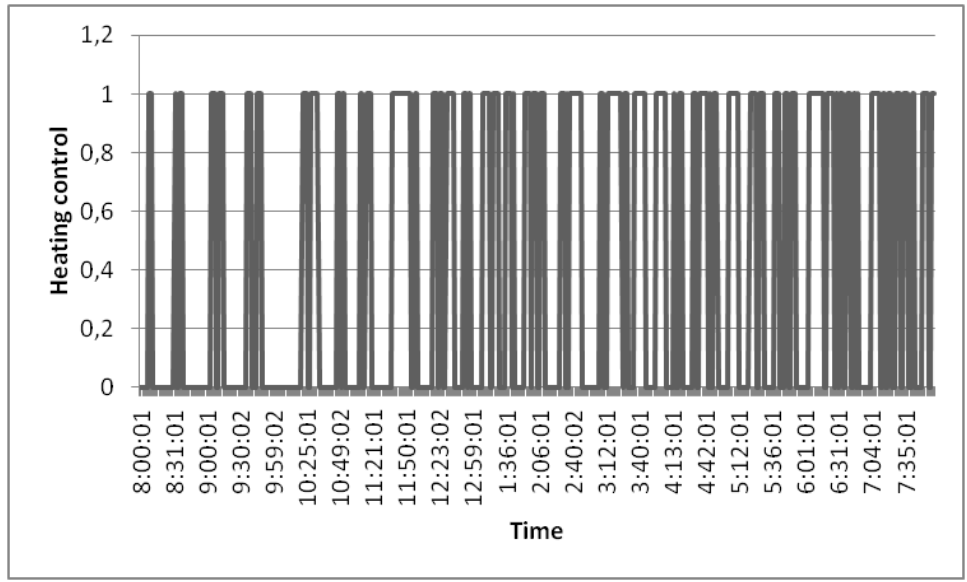}
      \caption{Heating control with model-free control}
      \label{Heating4}
   \end{figure}

An actuator fault can be described by Equation \eqref{acc1}.
An actuator fault on the heating control is simulated by a loss of efficiency equal to 50 \%. Figure \ref{Temperature4} shows results for the internal temperature with the temperature reference output equal to $18^{o}C$ during the night of 12-13 February 2014. Figure \ref{Heating4} demonstrates the accommodation ability of the heating control. The output temperature remains moreover very close of the internal temperature reference value.

\begin{figure}[!tbp]
      \centering
      \includegraphics[scale=0.6]{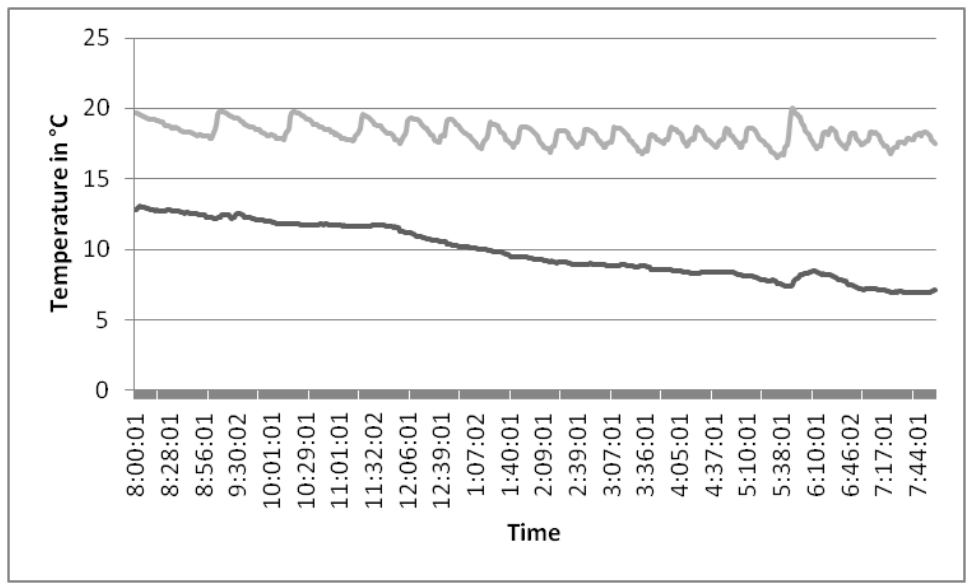}
      \caption{Temperature with model-free control (Te: Black line - Ti: Grey line)}
      \label{Temperature5}
   \end{figure}

\begin{figure}[!tbp]
      \centering
      \includegraphics[scale=0.6]{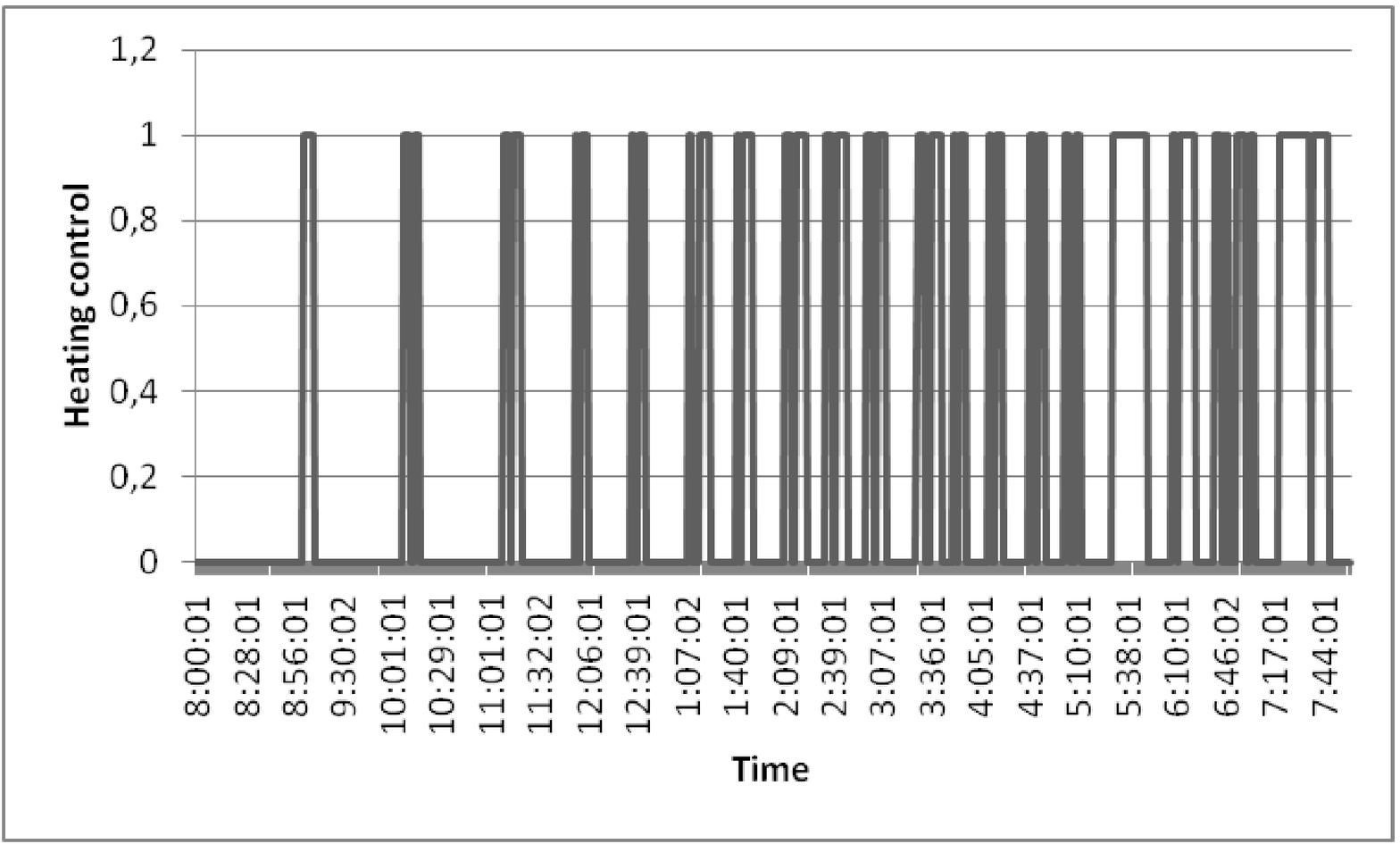}
      \caption{Heating control with model-free control}
      \label{Heating5}
   \end{figure}

Another actuator fault confirms the previous facts. Figure \ref{Temperature5} shows the results for the internal temperature with the temperature reference output equal to $18^{o}C$ during the night of 13-14 February 2014, with a loss of efficiency equal to 25 \%. The performances displayed by Figure \ref{Heating5} and Table \ref{Evaluation4} are again excellent.

\begin{table}[!tbp]
\caption{Results evaluation for the model-free control}
\label{Evaluation4}
\begin{center}
\begin{tabular}{|c||c||c|}
\hline
Output error & mean & variance\\
\hline
$e_{T_{i}}$ with $\beta$ = 50 \% & $-0.2^{o}C$ & $0.3^{o}C$\\
\hline
$e_{T_{i}}$ with $\beta$ = 25 \% & $0.2^{o}C$ & $0.5^{o}C$\\
\hline
\end{tabular}
\end{center}
\end{table}

\section{Conclusion} \label{conc}
Our successful model-free control strategy and its fault-tolerant capabilities will be further developed by taking advantage of technologically more advanced greenhouse systems.
{Let us mention here, among many other possibilities, a regulation of the $CO_{2}$ rate. Further comparisons with various other feedback synthesis techniques should also be investigated.}
We also hope that similar techniques might be useful in more or less analogous domains like air-conditioning in buildings (see, \emph{e.g.}, \cite{liu}). Data mining techniques will also be considered (see, \emph{e.g.}, \cite{dm}).


\section*{References}

\begin{thebibliography}{46}
\expandafter\ifx\csname natexlab\endcsname\relax\def\natexlab#1{#1}\fi
\providecommand{\url}[1]{\texttt{#1}}
\providecommand{\href}[2]{#2}
\providecommand{\path}[1]{#1}
\providecommand{\DOIprefix}{doi:}
\providecommand{\ArXivprefix}{arXiv:}
\providecommand{\URLprefix}{URL: }
\providecommand{\Pubmedprefix}{pmid:}
\providecommand{\doi}[1]{\href{http://dx.doi.org/#1}{\path{#1}}}
\providecommand{\Pubmed}[1]{\href{pmid:#1}{\path{#1}}}
\providecommand{\bibinfo}[2]{#2}
\ifx\xfnm\relax \def\xfnm[#1]{\unskip,\space#1}\fi
\bibitem[{Aaslyng \emph{et al.}(2005)Aaslyng, Ehler, Jakobsen}]{aa}
\bibinfo{author}{Aaslyng, J.M.}, \bibinfo{author}{Ehler, N.}, \bibinfo{author}{Jakobsen, L.},   \bibinfo{year}{2005}.
\newblock \bibinfo{title}{Climate control software integration with a greenhouse environmental control computer}
\newblock \bibinfo{journal}{Environ. Model. Soft.} \bibinfo{volume}{20},   \bibinfo{pages}{521--527}.
\bibitem[{Abouaissa \emph{et al.}(2012)Abouaissa, Fliess, Iordanova, Join}]{sofia}
\bibinfo{author}{Abouaissa, H.}, \bibinfo{author}{Fliess, M.},
  \bibinfo{author}{Iordanova, V.}, \bibinfo{author}{Join, C.},
  \bibinfo{year}{2012}.
\newblock \bibinfo{title}{Freeway ramp metering control made easy and
  efficient}.
\newblock \bibinfo{journal}{13$^{th}$ IFAC Symp. Control Transportation
  Systems, Sofia}.
\bibitem[{Arvantis  \emph{et al.}(2000)Arvantis, Paraskevopoulos and Vernados}]{c5}
\bibinfo{author}{Arvantis, K.}, \bibinfo{author}{Paraskevopoulos, P.},
  \bibinfo{author}{Vernados, A.}, \bibinfo{year}{2000}.
\newblock \bibinfo{title}{Multirate adaptative temperature control of
  greenhouses}.
\newblock \bibinfo{journal}{Comput. Electron Agric.} \bibinfo{volume}{26},
  \bibinfo{pages}{303--320}.
\bibitem[{{\AA}strom, H\"agglund(2006)}]{astrom}
\bibinfo{author}{{\AA}strom, K.}, \bibinfo{author}{H\"agglund, T.},
\bibinfo{year}{2006}.
\newblock \bibinfo{title}{Advanced PID Control}.
\newblock \bibinfo{journal}{Instrument Soc. Amer.}.
\bibitem[{Balmat, Lafont(2003)}]{c8}
\bibinfo{author}{Balmat, J.F.}, \bibinfo{author}{Lafont, F.},
  \bibinfo{year}{2003}.
\newblock \bibinfo{title}{Multi-model architecture supervised by Kohonen map}.
\newblock \bibinfo{journal}{Sci. Electron. Techno. Inform. Telecom (SETIT'03),
  Sousse} , \bibinfo{pages}{98--104}.
\bibitem[{Bennis \emph{et al.}(2008)Bennis, Duplaix, En\'ea, Haloua, Youlal}]{c7}
\bibinfo{author}{Bennis, N.}, \bibinfo{author}{Duplaix, J.},
  \bibinfo{author}{En\'ea, G.}, \bibinfo{author}{Haloua, M.},
  \bibinfo{author}{Youlal, H.}, \bibinfo{year}{2008}.
\newblock \bibinfo{title}{Greenhouse climate modelling and robust control}.
\newblock \bibinfo{journal}{Comput. Electron. Agric.}
  \bibinfo{volume}{61}, \bibinfo{pages}{96--107}.
\bibitem[{Blasco \emph{et al.}(2007)Blasco, Martinez, Herrero, Ramos,
  Sanchis}]{blasco}
\bibinfo{author}{Blasco, X.}, \bibinfo{author}{Martinez, M.},
  \bibinfo{author}{Herrero, J.}, \bibinfo{author}{Ramos, C.},
  \bibinfo{author}{Sanchis, J.}, \bibinfo{year}{2007}.
\newblock \bibinfo{title}{Model-based predictive control of greenhouse climate
  for reducing energy and water consumption}.
\newblock \bibinfo{journal}{Comput. Electron. Agric.}
  \bibinfo{volume}{55}, \bibinfo{pages}{49--70}.
\bibitem[{Bontsema \emph{et al.}(2011)Bontsema, Van Henten, Gieling,
  Swinkels}]{Bontsema}
\bibinfo{author}{Bontsema, J.}, \bibinfo{author}{Van Henten, E.},
  \bibinfo{author}{Gieling, T.}, \bibinfo{author}{Swinkels, G.},
  \bibinfo{year}{2011}.
\newblock \bibinfo{title}{The effect of sensor errors on production and energy
  consumption in greenhouse horticulture}.
\newblock \bibinfo{journal}{Comput. Electron. Agric.}
  \bibinfo{volume}{79}, \bibinfo{pages}{63--66}.
\bibitem[{Boulard(1989)}]{boul}
\bibinfo{author}{Boulard, T.}, \bibinfo{year}{1989}.
\newblock \bibinfo{title}{Water vapour transfer in a plastic house equipped
  with a dehumidification heat pump}.
\newblock \bibinfo{journal}{J. Agric. Engin. Res.} \bibinfo{volume}{44},
  \bibinfo{pages}{191--–204}.
\bibitem[{Callais(2006)}]{Callais}
\bibinfo{author}{Callais, M.J.}, \bibinfo{year}{2006}.
\newblock \bibinfo{title}{Les productions migrant entre terres et serres}.
\newblock \bibinfo{journal}{Agreste Primeur} \bibinfo{volume}{173}. Available at \newline {\tt http://agreste.agriculture.gouv.fr/IMG/pdf/primeur173.pdf}
 \bibitem[{Caponetto \emph{et al.}(2000)}]{cap}
\bibinfo{author}{Caponetto, R.}, \bibinfo{author}{Fortuna, L.}, \bibinfo{author}{Nunnari, G.}, \bibinfo{author}{Occhipinti, L.}, \bibinfo{author}{Xibilia, M.G.}, \bibinfo{year}{2000}.
\newblock \bibinfo{title}{Soft computing for greenhouse climate control}.
\newblock \bibinfo{journal}{IEEE Trans. Fuzzy Syst.} \bibinfo{volume}{8},
  \bibinfo{pages}{753--–760}.
\bibitem[{Cate, Challa(1984)}]{udink}
\bibinfo{author}{Cate, A.U.T.}, \bibinfo{author}{Challa, H.},
  \bibinfo{year}{1984}.
\newblock \bibinfo{title}{On optimal computer control of the crop growth
  system}.
\newblock \bibinfo{journal}{Acta Hortic.} \bibinfo{volume}{148},
  \bibinfo{pages}{267--276}.
\bibitem[{Critten, Bailey(2002)}]{critten}
\bibinfo{author}{Critten, D.}, \bibinfo{author}{Bailey, B.},
  \bibinfo{year}{2002}.
\newblock \bibinfo{title}{A review of greenhouse engineering developments
  during the 1990s}.
\newblock \bibinfo{journal}{Agric. Forest Meteorology} \bibinfo{volume}{112},
  \bibinfo{pages}{1--22}.
\bibitem[{Cunha \emph{et al.}(1997)Cunha, Couto, Ruano}]{c1}
\bibinfo{author}{Cunha, J.}, \bibinfo{author}{Couto, C.},
  \bibinfo{author}{Ruano, A.}, \bibinfo{year}{1997}.
\newblock \bibinfo{title}{Real-time parameter estimation of dynamic temperature
  models for greenhouse environmental control}.
\newblock \bibinfo{journal}{Control Eng. Practice} \bibinfo{volume}{5},
  \bibinfo{pages}{1473--1481}.
\bibitem[{De Miras \emph{et al.}(2013)De Miras, Join, Fliess, Riachy, Bonnet}]{compiegne}
\bibinfo{author}{De Miras, J.}, \bibinfo{author}{Join, C.},
  \bibinfo{author}{Fliess, M.}, \bibinfo{author}{Riachy, S.},
  \bibinfo{author}{Bonnet, S.}, \bibinfo{year}{2013}.
\newblock \bibinfo{title}{Active magnetic bearing: A new step for model-free control}.
\newblock \bibinfo{journal}{52$^{nd}$ IEEE Conf. Decision Control, Florence}. Preprint available at \newline
{\tt http://hal.archives-ouvertes.fr/hal-00857649/en/}
\bibitem[{Dong \emph{et al.}(2013)}]{dong}
\bibinfo{author}{Dong, Q.}, \bibinfo{author}{Yang, W.},
 \bibinfo{author}{Yang, L.}, \bibinfo{author}{Chen, S.},
  \bibinfo{author}{Du, S.}, \bibinfo{author}{Li, F.}, \bibinfo{author}{Shi, Q.}, \bibinfo{author}{Xu, Y.}, \bibinfo{year}{2013}.
\newblock \bibinfo{title}{Crop model-based greenhouse optimal control system: Survey and perspectives}.
\newblock \bibinfo{journal}{D. Li, Y. Chen (Eds): Computer  Computing Techno. Agric. VI --
IFIP Advan. Informat. Communicat. Techno.} \bibinfo{volume}{392},  \bibinfo{pages}{216-224},
\newblock \bibinfo{journal}{Springer}.
\bibitem[{Duarte-Galvan \emph{et al.}(2012)Duarte-Galvan, Torres-Pacheco,
  Guevara-Gonzalez, Romero-Troncoso, Contreras-Medina, Rios-Alcaraz and
  Millan-Almaraz}]{c19}
\bibinfo{author}{Duarte-Galvan, C.}, \bibinfo{author}{Torres-Pacheco, I.},
  \bibinfo{author}{Guevara-Gonzalez, R.}, \bibinfo{author}{Romero-Troncoso, R.},
  \bibinfo{author}{Contreras-Medina, L.}, \bibinfo{author}{Rios-Alcaraz, M.},
  \bibinfo{author}{Millan-Almaraz, J.}, \bibinfo{year}{2012}.
\newblock \bibinfo{title}{Advantages and disadvantages of control theories
  applied in greenhouse climate control systems}.
\newblock \bibinfo{journal}{Spanish J. Agri. Res.} \bibinfo{volume}{10},
  \bibinfo{pages}{926--938}.
\bibitem[{El Ghoumari \emph{et al.}(2005)El Ghoumari, Tantau,
 Serrano}]{el}
\bibinfo{author}{El Ghoumari, M.Y.}, \bibinfo{author}{Tantau, H.-J.},
  \bibinfo{author}{Serrano, J.}, \bibinfo{year}{2005}.
\newblock \bibinfo{title}{
Non-linear constrained MPC: Real-time implementation of greenhouse air temperature control}.
\newblock \bibinfo{journal}{Comput. Electron. Agricul.} \bibinfo{volume}{49},
  \bibinfo{pages}{345--356}.
\bibitem[{Fliess, Join(2013)Fliess and Join}]{ijc13}
\bibinfo{author}{Fliess, M.}, \bibinfo{author}{Join, C.}, \bibinfo{year}{2013}.
\newblock\bibinfo{title} Model-free control,
\newblock\bibinfo{title}{Int. J. Contr.}  \bibinfo{volume}86, \bibinfo{pages}{2228--2252}.
\bibitem[{Fliess, Join(2014)Fliess and Join}]{ecc}
\bibinfo{author}{Fliess, M.}, \bibinfo{author}{Join, C.}, \bibinfo{year}{2014}.
\newblock\bibinfo{title}Stability margins and model-free control: A first look,
\newblock \bibinfo{journal}{13$^{th}$ Europ. Contr. Conf., Strasbourg}. Preprint available at  \newline
{\tt http://hal.archives-ouvertes.fr/hal-00966107/en/}
\bibitem[{Fliess \emph{et al.}(2004)Fliess, Join, Sira-Ram\'{\i}rez}]{ijc04}
\bibinfo{author}{Fliess, M.}, \bibinfo{author}{Join, C.},
\bibinfo{author}{Sira-Ram\'{\i}rez, H.}, \bibinfo{year}{2004}.
\newblock \bibinfo{title}{Robust residual generation for linear fault diagnosis: an algebraic setting with examples}.
\newblock \bibinfo{journal}{Int. J. Contr.} \bibinfo{year}{77}, \bibinfo{pages}{1223--1242}.
\bibitem[{Fliess \emph{et al.}(2008)Fliess, Join, Sira-Ram\'{\i}rez}]{nl}
\bibinfo{author}{Fliess, M.}, \bibinfo{author}{Join, C.},
\bibinfo{author}{Sira-Ram\'{\i}rez, H.}, \bibinfo{year}{2008}.
\newblock \bibinfo{title}{Non-linear estimation is easy}.
\newblock \bibinfo{journal}{Int. J. Model. Identif. Control}
\bibinfo{volume}{4}, \bibinfo{pages}{12--27}.
\bibitem[{Fliess, Sira-Ram\'{\i}rez(2003)}]{sira1}
\bibinfo{author}{Fliess, M.}, \bibinfo{author}{Sira-Ram\'{\i}rez, H.},
  \bibinfo{year}{2003}.
\newblock \bibinfo{title}{An algebraic framework for linear identification}.
\newblock \bibinfo{journal}{ESAIM Control Optimiz. Calc. Variat.}
  \bibinfo{volume}{9}, \bibinfo{pages}{151--168}.
\bibitem[{Fliess, Sira-Ram\'{\i}rez(2008)}]{sira2}
\bibinfo{author}{Fliess, M.}, \bibinfo{author}{Sira-Ram\'{\i}rez, H.},
  \bibinfo{year}{2008}.
\newblock \bibinfo{title}{Closed-loop parametric identification for
  continuous-time linear systems via new algebraic techniques}.
\newblock \bibinfo{journal}{Eds. H. Garnier and L. Wang: Identification of
  Continuous-time Models from Sampled Data, Springer},
  \bibinfo{pages}{362--391}.
\bibitem[{Fourati(2014)Fourati}]{fourati}
\bibinfo{author}{Fourati, F.},
\bibinfo{year}{2014}.
\newblock \bibinfo{title}{Multiple neural control of a greenhouse}.
\newblock \bibinfo{journal}{Neurocomput.} \bibinfo{volume}{139},
\bibinfo{pages}{138--144}.
\bibitem[{Gao(2014)}]{gao}
\bibinfo{author}{Gao, Z.},
\bibinfo{year}{2014}.
\newblock \bibinfo{title}{On the centrality of disturbance rejection in automatic control}.
\newblock \bibinfo{journal}{ISA Trans.} \bibinfo{volume}{53}, \bibinfo{pages}{850-857}.
\bibitem[{G\'edouin \emph{et al.}(2011)G\'{e}douin, Delaleau, Bourgeot, Join, Arab-Chirani, Calloch}]{brest}
\bibinfo{author}{G\'{e}douin, P.-A.}, \bibinfo{author}{Delaleau, E.}, \bibinfo{author}{Bourgeot, J.-M.}, \bibinfo{author}{Join, C.},
\bibinfo{author}{Arab-Chirani, S.}, \bibinfo{author}{Calloch, S.}, \bibinfo{year}{2011}.
\newblock \bibinfo{title}{Experimental comparison
of classical pid and model-free control: position control of a shape
memory alloy active spring}.
\newblock \bibinfo{journal}{Control Eng. Practice} \bibinfo{volume}{19},
  \bibinfo{pages}{433--441}.
\bibitem[{Gruber \emph{et al.} (2011)Gruber, Guzm\'an, Rodr\'{\i}guez, Bordons,
  Berenguel, S\'anchez}]{gruber}
\bibinfo{author}{Gruber, J.}, \bibinfo{author}{Guzm\'an, J.},
  \bibinfo{author}{Rodr\'{\i}guez, F.}, \bibinfo{author}{Bordons, C.},
  \bibinfo{author}{Berenguel, M.}, \bibinfo{author}{S\'anchez, J.},
  \bibinfo{year}{2011}.
\newblock \bibinfo{title}{Nonlinear MPC based on a Volterra series model for
  greenhouse temperature control using natural ventilation}.
\newblock \bibinfo{journal}{Control Eng. Practice} \bibinfo{volume}{19},
  \bibinfo{pages}{354--366}.
\bibitem[{Hou \emph{et al.}(2006)Hou, Lian, Yuan}]{dm}
\bibinfo{author}{Hou, Z.}, \bibinfo{author}{Lian, Z.}, \bibinfo{author}{Yao,
  Y.}, \bibinfo{author}{Yuan, X.}, \bibinfo{year}{2006}.
\newblock \bibinfo{title}{Data mining based sensor fault diagnosis and
  validation for building air conditioning system}.
\newblock \bibinfo{journal}{Energy Convers. Manag.} \bibinfo{volume}{47},
  \bibinfo{pages}{2479--2490}.
\bibitem[{Ioslovich \emph{et al.}(2009)Ioslovich, Gutman, Linker}]{io}
\bibinfo{author}{Ioslovich, I.}, \bibinfo{author}{Gutman, P.O.}, \bibinfo{author}{Linker,
  R.}, \bibinfo{year}{2009}.
\newblock \bibinfo{title}{Hamilton-Jacobi-Bellman formalism for optimal climate control of
greenhouse crop}.
\newblock \bibinfo{journal}{Automatica} \bibinfo{volume}{45},
  \bibinfo{pages}{1227-1231}.
\bibitem[{Isermann(2011)}]{isermann}
\bibinfo{author}{Isermann, R.}, \bibinfo{year}{2011}.
\newblock \bibinfo{title}{Fault-Diagnosis Applications}.
\newblock \bibinfo{journal}{Springer}.
\bibitem[{Join \emph{et al.}(2013)Join, Robert and Fliess}]{nice}
\bibinfo{author}{Join, C.}, \bibinfo{author}{Chaxel, F.}, \bibinfo{author}{Fliess, M.}, \bibinfo{year}{2013}.
\newblock \bibinfo{title}{``Intelligent'' controllers on cheap and small programmable devices}.
\newblock \bibinfo{journal}{2$^{nd}$ Int. Conf. Control Fault-Tolerant Syst., Nice}. Preprint available at \newline
{\tt http://hal.archives-ouvertes.fr/hal-00845795/en/}
\bibitem[{Kiltz \emph{et al.}(2014)Kiltz, Join, Mboup, Rudolph}]{sarre}\bibinfo{author}{Kiltz, L.}, \bibinfo{author}{Join, C.}, \bibinfo{author}{Mboup, M.}, \bibinfo{author}{Rudolph, J.},  \bibinfo{year}{2014}.
\newblock \bibinfo{title}{Fault-tolerant control based on algebraic derivative estimation applied on a magnetically supported plate}. \newblock \bibinfo{journal}{Contr. Engin. Practice}  \bibinfo{volume}{26}, \bibinfo{pages}{107--115}.
\bibitem[{Kimball(1973)}]{c16}
\bibinfo{author}{Kimball, B.}, \bibinfo{year}{1973}.
\newblock \bibinfo{title}{Simulation of the energy balance of a greenhouse}.
\newblock \bibinfo{journal}{Agric. Meteorology} \bibinfo{volume}{11},
  \bibinfo{pages}{243--260}.
\bibitem[{Kittas, Batzanas(2010)}]{kit}
\bibinfo{author}{Kittas, C.}, \bibinfo{author}{Batzanas, T.},
\bibinfo{year}{2007}.
\newblock \bibinfo{title}{Greenhouse microclimate and dehumidification effectiveness under different ventilator configurations}.
\newblock \bibinfo{journal}{Building Environ.} \bibinfo{volume}{42}, \bibinfo{pages}{3774--3784}
\bibitem[{Lafont, Balmat(2002)}]{c6}
\bibinfo{author}{Lafont, F.}, \bibinfo{author}{Balmat, J.F.},
  \bibinfo{year}{2002}.
\newblock \bibinfo{title}{Optimized fuzzy control of a greenhouse}.
\newblock \bibinfo{journal}{Fuzzy Sets Syst.} \bibinfo{volume}{128},
  \bibinfo{pages}{47--59}.
\bibitem[{Lafont \emph{et al.}(2005)Lafont, Balmat and Taurines}]{c3}
\bibinfo{author}{Lafont, F.}, \bibinfo{author}{Balmat, J.F.},
  \bibinfo{author}{Taurines, M.}, \bibinfo{year}{2005}.
\newblock \bibinfo{title}{Fuzzy forgetting factor for system identification}.
\newblock \bibinfo{journal}{3$^{rd}$ IEEE Int. Conf. Systems Signals Devices,
  Sousse, 2005}.
\bibitem[{Lafont \emph{et al.}(2014)Lafont, Pessel, Fliess}]{Laflast2}
\bibinfo{author}{Lafont, F.}, \bibinfo{author}{Balmat, J.F.},
  \bibinfo{author}{Pessel, N.}, \bibinfo{author}{Fliess, M.},
  \bibinfo{year}{2014}.
\newblock \bibinfo{title}{Model-free control and fault accommodation for an
  experimental greenhouse}.
\newblock \bibinfo{journal}{International Conference on Green Energy and
  Environmental Engineering, Sousse}. Preprint available at
 {\tt http://hal.archives-ouvertes.fr/hal-00978226/en/}
\bibitem[{Lafont \emph{et al.}(2013)Lafont, Pessel, Balmat and Fliess}]{Laflast}
\bibinfo{author}{Lafont, F.}, \bibinfo{author}{Pessel, N.},
  \bibinfo{author}{Balmat, J.F.}, \bibinfo{author}{Fliess, M.},
  \bibinfo{year}{2013}.
\newblock \bibinfo{title}{On the model-free control of an experimental
  greenhouse}.
\newblock \bibinfo{journal}{World Congress Engineering Computer Science
  2013, International Association of Engineers, San Francisco}. Preprint available at
 {\tt http://hal.archives-ouvertes.fr/hal-00831598/en/}
 \bibitem[{Landau \emph{et al.}(2011)Landau, Lozano, M'Saad, Karimi}]{landau}
\bibinfo{author}{Landau, I.D.}, \bibinfo{author}{Lozano, R.}, \bibinfo{author}{M'Saad, M.}, \bibinfo{author}{Karimi, A.}, \bibinfo{year}{2011}.
\newblock \bibinfo{title}{Adaptive Control}.
\newblock \bibinfo{journal}{2$^{nd}$ ed., Springer}.
\bibitem[{Liu \emph{et al.}(2013)Liu, Jiang, Zhang}]{liu} \bibinfo{author}{Liu X.}, \bibinfo{author}{Jiang Y.}, \bibinfo{author}{Zhang T.},  \bibinfo{year}{2013}.
\newblock \bibinfo{title}{Temperature and Humidity Independent Control (THIC) of Air-conditionning System}.
\newblock \bibinfo{journal}{Springer}.


\bibitem[{Medjber(2012)}]{eue}
\bibinfo{author}{Medjber, A.}, \bibinfo{year}{2012}.
\newblock \bibinfo{title}{Automatisation d'une serre agricole: Commande et r\'{e}gulation}.
\newblock \bibinfo{journal}{\'Edit. Univ. Europ.}.

\bibitem[{Menhour \emph{et al.}(2013)}]{menhour}
\bibinfo{auteur}{Menhour, L.}, \bibinfo{auteur}{D'Andr\'ea-Novel, B.}, \bibinfo{auteur}{Fliess, M.}, \bibinfo{auteur}{Mounier, H.}, \bibinfo{year}{2013}.
\newblock\bibinfo{title}{Multivariable decoupled longitudinal and lateral vehicle control: A model-free design}.
\newblock\bibinfo{journal}{52$^{nd}$ IEEE Conf. Decision Control, Florence}. Preprint available at \newline
{\tt http://hal.archives-ouvertes.fr/hal-00859444/en/}

\bibitem[{Noura \emph{et al.}(2009)Noura,Theilliol,Ponsart,Chamseddine}]{nancy}
\bibinfo{author}{Noura, H.}, \bibinfo{author}{Theilliol, D.}, \bibinfo{author}{Ponsart, J.-C.}, \bibinfo{author}{Chamseddine, A.}, \bibinfo{year}{2009}.
\newblock \bibinfo{title}{Fault-tolerant Control Systems: Design and Practical Applications}.
\bibinfo{journal}{Springer}.
\bibitem[{O'Dwyer(2009)}]{od}
\bibinfo{author}{O'Dwyer, A.}, \bibinfo{year}{2009}.
\newblock \bibinfo{title}{Handbook of PI and PID Controller Tuning Rules}.
\newblock \bibinfo{journal}{(3$^{rd}$ ed.), London : Imperial College Press}.
\bibitem[{Pasgianos \emph{et al.}(2003)Pasgianos, Arvanitis, Polycarpou, Sigrimis}]{pas}
\bibinfo{author}{Pasgianos,  G.D.}, \bibinfo{author}{Arvanitis, K.G.}, \bibinfo{author}{Polycarpou, P.}, \bibinfo{author}{Sigrimis, N.}, \bibinfo{year}{2003}.
\newblock \bibinfo{title}{A nonlinear feedback technique for greenhouse environmental control}.
\newblock \bibinfo{journal}{Comput. Electron. Agricul.} \bibinfo{volume}{40},
\bibinfo{pages}{153--177}.
\bibitem[{Pessel, Balmat(2005)}]{c18}
\bibinfo{author}{Pessel, N.}, \bibinfo{author}{Balmat, J.F.},
  \bibinfo{year}{2005}.
\newblock \bibinfo{title}{Principal component analysis to the modelling of
  systems - application to an experimental greenhouse}.
\newblock \bibinfo{journal}{3$^{rd}$ IEEE Int. Conf. Systems Signals Devices,
  Sousse}.
\bibitem[{Pessel \emph{et al.}(2009)Pessel, Duplaix, Balmat and Lafont}]{c4}
\bibinfo{author}{Pessel, N.}, \bibinfo{author}{Duplaix, J.},
  \bibinfo{author}{Balmat, J.F.}, \bibinfo{author}{Lafont, F.},
  \bibinfo{year}{2009}.
\newblock \bibinfo{title}{A multi-structure modeling methodology}.
\newblock \bibinfo{journal}{V.E. Balas, J. Fodor, A.R.
  V\'{a}rkonyi-K\'{o}czy, Eds: Soft Computing Based Modeling in Intelligent Systems,
  Studies in Computational Intelligence, Springer}, \bibinfo{volume}{196},
  \bibinfo{pages}{93--113}.
\bibitem[{Pi\~n\'on \emph{et al.}(2005)Pi\~n\'on, Camacho, Kuchen, Pe\~na}]{Pin}
\bibinfo{author}{Pi\~n\'on, S.}, \bibinfo{author}{Camacho, E.F.},
  \bibinfo{author}{Kuchen, B.}, \bibinfo{author}{Pe\~na, M.},
  \bibinfo{year}{2005}.
\newblock \bibinfo{title}{Constrained predictive control of a greenhouse}.
\newblock \bibinfo{journal}{Comput. Electron. Agricul.} \bibinfo{volume}{49},
  \bibinfo{pages}{317--329}.
  \bibitem[{Ponce \emph{et al.}(2012)}]{crc14}
\bibinfo{author}{Ponce, P.}, \bibinfo{author}{Molina, A.},
 \bibinfo{author}{Cepeda, P.}, \bibinfo{author}{Lugo, E.},
 \bibinfo{year}{2014}.
\newblock \bibinfo{title}{Greenhouse Design and Control}.
\newblock \bibinfo{journal}{CRC Press}.
\bibitem[{Ram\'{\i}rez-Arias \emph{et al.}(2012)Ram\'{\i}rez-Arias, Rodr\'{\i}guez, Guzm\'an, Berenguel}]{automatica}
\bibinfo{author}{Ram\'{\i}rez-Arias, J.L.}, \bibinfo{author}{Rodr\'{\i}guez, F.},
 \bibinfo{author}{Guzm\'an, J.}, \bibinfo{author}{Berenguel, M.},
 \bibinfo{year}{2012}.
\newblock \bibinfo{title}{Multiobjective hierarchical control architecture for greenhouse crop growth}.
\newblock \bibinfo{journal}{Automatica} \bibinfo{volume}{48},
  \bibinfo{pages}{490--498}.
\bibitem[{Rodr\'{\i}guez \emph{et al.}(2015)Rodr\'{\i}guez, Berenguel, Guzman, Ram\'{\i}rez-Arias}]{rod}
\bibinfo{author}{Rodr\'{\i}guez, F.}, \bibinfo{author}{Berenguel, M.}, \bibinfo{author}{Guzman, J.L.}, \bibinfo{author}{Ram\'{\i}rez-Arias, A.},   \bibinfo{year}{2015}.
\newblock \bibinfo{title}{Modeling and Control of Greenhouse Crop Growth}.
\newblock \bibinfo{journal}{Springer}.
\bibitem[{Salgado, Cunha(2005)}]{c2}
\bibinfo{author}{Salgado, P.}, \bibinfo{author}{Cunha, J.},
  \bibinfo{year}{2005}.
\newblock \bibinfo{title}{Greenhouse climate hierarchical fuzzy modelling}.
\newblock \bibinfo{journal}{Control Eng. Practice} \bibinfo{volume}{13},
  \bibinfo{pages}{613--628}.
\bibitem[{Shamshiri, Wan Ismail(2013)}]{sham}
\bibinfo{author}{Shamshiri, R.}, \bibinfo{author}{Ismail, W.I.W.},
  \bibinfo{year}{2013}.
\newblock \bibinfo{title}{A review of greenhouse climate control and automation systems in tropical regions}.
\newblock \bibinfo{journal}{J. Agric. Sci. Appl.} \bibinfo{volume}{2},
  \bibinfo{pages}{176--183}.
\bibitem[{Shumsky \emph{et al.}(2011)Shumsky, Zhirabok, Jiang}]{shum}
\bibinfo{author}{Shumsky, A.}, \bibinfo{author}{Zhirabok, A.},
  \bibinfo{author}{Jiang, B.}, \bibinfo{year}{2011}.
\newblock \bibinfo{title}{Fault accomodation in nonlinear and linear dynamic
  systems: Fault decoupling based approach}.
\newblock \bibinfo{journal}{Int. J. Innovat. Comput. Informat. Contr.}
  \bibinfo{volume}{7}, \bibinfo{pages}{4535--4549}.
 \bibitem[{ Sira-Ram\'{\i}rez \emph{et al.}(2014)Sira-Ram\'{\i}rez, Garc\'{\i}a-Rodr\'{\i}guez, Cort\`{e}s-Romero, Luviano-Ju\'{a}rez}]{sira}
 \bibinfo{author}{Sira-Ram\'{\i}rez, H.}, \bibinfo{author}{Garc\'{\i}a-Rodr\'{\i}guez, C.}, \bibinfo{author}{Cort\`{e}s-Romero, J.}, \bibinfo{author}{Luviano-Ju\'{a}rez,  A.}, \bibinfo{year}{2014},
 \newblock \bibinfo{title}{Algebraic Identification and Estimation Methods in Feedback Control Systems}.
 \newblock \bibinfo{journal}{Wiley}.
  \bibitem[{ Speetjens \emph{et al.}(2009)Speetjens, Stigter, Van Straten}]{speetens}
 \bibinfo{author}{Speetjens, S.L.},  \bibinfo{author}{Stigter, J.D.}, \bibinfo{author}{Van Straten, G.}, \bibinfo{year}{2009},
 \newblock \bibinfo{title}{Towards an adaptive model for greenhouse control}.
 \newblock \bibinfo{journal}{Comput. Electron. Agricul.} \bibinfo{volume}{67},
  \bibinfo{pages}{1--8}.
\bibitem[{van Straten \emph{et al.}(2010)van Straten, van Willigenburg, van Henten, van Ooteghem}]{tf}
\bibinfo{author}{van Straten, G.}, \bibinfo{author}{van Willigenburg, H.},
  \bibinfo{author}{van Henten, F.},   \bibinfo{author}{van Ooteghem, R.}, \bibinfo{year}{2010}.
\newblock \bibinfo{title}{Optimal Control of Greenhouse Cultivation}.
\newblock \bibinfo{journal}{CRC Press}.
\bibitem[{Tchamitchian \emph{et al.}(2006)Tchamitchian, Martin-Clouaire, Lagier,
  Jeannequin and Mercier}]{tcham}
\bibinfo{author}{Tchamitchian, M.}, \bibinfo{author}{Martin-Clouaire, R.},
  \bibinfo{author}{Lagier, J.}, \bibinfo{author}{Jeannequin, J.},
  \bibinfo{author}{Mercier, S.}, \bibinfo{year}{2006}.
\newblock \bibinfo{title}{Serriste: A daily set point determination software
  for glasshouse tomato production}.
\newblock \bibinfo{journal}{Comput. Electron. Agric.}
  \bibinfo{volume}{50}, \bibinfo{pages}{25--47}.
\bibitem[{Urban \emph{et al.}(2010)}]{c15}
\bibinfo{author}{Urban, L.}, \bibinfo{author}{Urban, I.}, \bibinfo{year}{2010}.
\newblock \bibinfo{title}{La production sous serre, tome 1, la gestion du
  climat}.
\newblock \bibinfo{journal}{2$^{e}$ \'ed., Lavoisier}.
\bibitem[{Viard-Gaudin(1981)}]{c17}
\bibinfo{author}{Viard-Gaudin, C.}, \bibinfo{year}{1981}.
\newblock \bibinfo{title}{Simulation et commande auto-adaptative d'une serre
  agricole}.
\newblock \bibinfo{journal}{Th\`ese, Universit\'e de Nantes}.
\bibitem[{Villagra \emph{et al.}(2011a)Villagra, d'Andr\'ea-Novel, Fliess, Mounier}]{mines1}
\bibinfo{author}{Villagra, J.}, \bibinfo{author}{d'Andr\'{e}a-Novel, B.},
\bibinfo{author}{Fliess, M.}, \bibinfo{author}{Mounier, M.},
\bibinfo{year}{2011a}.
\newblock \bibinfo{title}{A diagnosis-based approach for tire-road forces and maximum friction estimation}.
\newblock \bibinfo{journal}{Contr. Engin. Practice}
\bibinfo{volume}{19}, \bibinfo{pages}{174--184}.
\bibitem[{Villagra \emph{et al.}(2011b)Villagra, d'Andr\'ea-Novel, Fliess, Mounier, Menhour}]{mines2}
\bibinfo{author}{Villagra, J.}, \bibinfo{author}{d'Andr\'{e}a-Novel, B.},
\bibinfo{author}{Fliess, M.}, \bibinfo{author}{Mounier, M.}, \bibinfo{author}{Menhour, L.},
\bibinfo{year}{2011b}.
\newblock \bibinfo{title}{Corrigendum to ``A diagnosis-based approach for tire-road forces and maximum friction estimation''}.
\newblock \bibinfo{journal}{Contr. Engin. Practice}
\bibinfo{volume}{19}, \bibinfo{pages}{1252}.
\bibitem[{Von Zabeltitz(2011)}]{von}
\bibinfo{author}{Von Zabeltitz, C.}, \bibinfo{year}{2011}.
\newblock \bibinfo{title}{Integrated Greenhouse Systems for Mild Climates}.
\newblock \bibinfo{journal}{Springer}.
\bibitem[{Yosida(1984)}]{yosida}
\bibinfo{author}{Yosida, K.}, \bibinfo{year}{1984}.
\newblock \bibinfo{title}{Operational Calculus (translated from the Japanese)}.
\newblock \bibinfo{journal}{Springer}.
\bibitem[{Zhang(2008)}]{zh}
\bibinfo{author}{Zhang, Z.}, \bibinfo{year}{2008}.
\newblock \bibinfo{title}{Multiobjective optimization immune algorithm in dynamic environments and its application to greenhouse control}.
\newblock \bibinfo{journal}{Appl. Soft Comput.} \bibinfo{volume}{8}, \bibinfo{pages}{959--971}.
\end{thebibliography}

\end{document}